\documentclass[journal]{new-aiaa}
\usepackage[utf8]{inputenc}

\usepackage{graphicx}
\usepackage{amsmath}
\usepackage[version=4]{mhchem}
\usepackage{siunitx}
\usepackage{longtable,tabularx}
\usepackage{pstricks}

\hypersetup{bookmarksnumbered, pdfstartview={FitH}, urlcolor=cyan,
 citecolor=blue, menucolor=blue, pagecolor=blue, linkcolor=blue}

\usepackage{amsmath,amssymb}
\usepackage{amsfonts,verbatim}
\usepackage{threeparttable,booktabs,dcolumn}
\usepackage{dsfont}

\usepackage{courier}
\usepackage[T1]{fontenc}

 \newtheorem{problem}{Problem}[section]

 \newtheorem{exam}{Example}
 
 \newtheorem{conclusion}{Claim}
 \newtheorem{rem}{Remark}
 \newtheorem{conj}{Conjecture}
 \newtheorem{bcs}{\rm List\rm}

\newcommand{\diag}{\mathop{\rm diag}}
 \newcommand{\abs}[1]{\left\vert#1\right\vert}

\usepackage{multirow,booktabs}
\usepackage{bm}

\usepackage{commath}
\renewcommand\thesection{\arabic{section}}

\title{Optimal Two-impulse Space Interception with Multi-constraints\footnote{The first version was 
completed in Dec.~2017.}}

\author{Li Xie\footnote{Professor,
School of Control and Computer Engineering; the correspondence author, lixie@ncepu.edu.cn
}}
\affil{State Key Laboratory of Alternate Electrical Power System with Renewable Energy Sources\\North China Electric Power University, Beijing 102206, P.R.~China}
\author{Yiqun Zhang\footnote{Senior Research Scientist, yiqunzhang@hotmail.com} and Junyan Xu\footnote{Associate Research Scientist, junyan\_Xu@sina.cn}
}
\affil{Beijing Institute of Electronic Systems Engineering, Beijing 100854, P.R.~China}

\begin{document}

\maketitle

\begin{abstract}

We consider optimal two-impulse space interception problems with multi-constraints. 
The multi-constraints are imposed on the terminal position of an interceptor, impulse and impact instants, and the component-wise magnitudes of velocity impulses.
We formulate these optimization problems as multi-point boundary value problems and the calculus of variations is used to solve them. All inequality constraints are converted into equality constraints by using slackness variable methods in order to use Lagrange multiplier method. 
A new dynamic slackness variable method is presented.  As a result, an  indirect optimization method is established for two-impulse space interception problems with multi-constraints. 
Subsequently, our method is used to solve the  two-impulse space interception problems of free-flight ballistic missiles. 
A number of conclusions have been established based on highly accurate numerical solutions.  Specifically, by numerical examples, we show that when time and velocity impulse constraints are imposed, optimal two-impulse  solutions may occur, and also if two impulse instants are free, then two-impulse space interception problems with velocity impulse constraints may degenerate to the one-impulse case. 
\end{abstract}



\section*{Nomenclature}


{\renewcommand\arraystretch{1.0}
\noindent\begin{longtable*}{@{}l @{\quad=\quad} l@{}}
$\abs \cdot$, $I_{3\times 3}$, BCs, diag  & magnitude,  the $3\times 3$ identical matrix, boundary conditions, and diagonal matrix \\
\multicolumn{2}{@{}l}{Vectors in boldface}\\
$\mathbf p$ &  costate or Lagrange multiplier in \eqref{eq_Dec6_01}\\
$ -\mathbf p_{Mv} $ & primer vector in \eqref{eq_Dec6_01}\\
$\mathbf q$ & Lagrange multiplier in \eqref{eq_Dec6_01}\\

$\mathbf r $  & position vector,  m\\
$\mathbf r_f$  & reference position vector, m \\
$\mathbf v, \Delta \mathbf v$  & velocity vector,  velocity impulse, m/s\\


$\mathbf x(t)$ & state, consists of $\mathbf r(t)$ and $\mathbf v(t)$\\
\multicolumn{2}{@{}l}{Scales and Greek Symbols}\\
$H, \hat H, \Delta \hat H$  & Hamiltonian functions  in \eqref{eq_7_2_1_Sep19}, \eqref{eq_Dec27_02}, \eqref{eq_Dec12_01},\\
$J$  & cost functional in \eqref{eq_Nov14_01}\\
$\tilde J, \hat J$  & augmented cost functionals in \eqref{eq_Dec6_01} and \eqref{eq_Dec27_01} or \eqref{eq_Dec15_01}\\
$k_3, k_4$ & weighting coefficients in \eqref{eq_Dec15_01}\\
$t_0, t_f$ & initial time, terminal time, s\\
$t_i, t_h$ & i-{\em th} impulse instant, impact instant, s \\

$t^-, t^+$ & time just before $t$,  time just after $t$\\
$ p_{i\min}, p_{i\max}$ & left and right boundary points of an interval in velocity   impulse constraints \eqref{eq_Nov5_01o}\\
$p_\epsilon, p_{i\epsilon}$  & costate of $\epsilon$,  costate of $\epsilon$ in the $i$-{\em th} time interval in \eqref{eq_May6_01}\\
$ r_{x\min}, r_{x\max}$ & lower or upper bounds in final position constraints \eqref{eq_Nov5_01m} in the direction of $x$\\

$R_e$  &  earth's radius 6378145 m in Appendix B\\
$(H, i, \Omega, e, \omega)$ & orbital elements in Appendix B\\
& altitude, inclination, right ascension of ascending node, eccentricity, and argument of perigee\\
$\alpha, \beta, \gamma, \eta$ & upper or lower bounds of time constraints in \eqref{eq_Nov5_01n}\\
$\delta$ & the first variation in Appendix A\\
$\epsilon$ & static or dynamic slackness variable in Section \ref{sect_3}\\
$\gamma_h$ & Lagrange multiplier w.r.t.~the interception condition in \eqref{eq_Dec6_01} \\
$\lambda, \lambda_i, \mu_i $ & Lagrange multipliers in \eqref{eq_Dec27_01},  \eqref{eq_Dec15_01} \\
$\mu$ & gravitational constant 3.986e+14\\
$\theta$ & true anomaly in Appendix B\\
\multicolumn{2}{@{}l}{Subscripts}\\
$*$ & optimal value\\
$0, f$ & initial, final or reference point\\
$i$ & index\\
$M, T$ & interceptor, target\\
$r, v, \epsilon$ & position, velocity, and slackness variable \\
$x,y,z$ & components of a vector\\
\multicolumn{2}{@{}l}{Superscripts}\\
$\text{T}$ & transposition of a matrix \\
\multicolumn{2}{@{}l}{In verbatim environment}\\
{\texttt{abs}}  & magnitude \\
{\texttt{dv1,dv2}}  &  velocity impulse vectors\\
{\texttt{t1,t2,th}}  & scaled time instants \\
{\texttt{t\_impulse}, \texttt{t\_impact}, \texttt{tf}}  & impulse instant, impact instant, and terminal instant \\
\end{longtable*}}

\section{Introduction} \label{Sec_0}

It is well-known that from the first order necessary  optimality, the calculus of variations can
transform an optimization problem into a two-point or multi-point boundary value problem.  
Arthur E. Bryson Jr.~\cite{Bryson_story} told such a story about how the calculus of variations was involved in aerospace field through the maximum range problem of a Hughes air-to-air missile 
in 1952. 
Using the variational method to solve aerospace problems can be traced back to the rocket age \cite{Bryson1996}.  For example, Hamel in 1927 formulated the well-known Goddard problem, i.e., optimizing the altitude of a rocket given the amount of propellant,  as a variational problem, and then Tsien and Evans \cite{Tsien} in 1951 re-considered Hamel's variational problem and gave an important analytical solution and also numerical data for two kinds of aerodynamic drags.
A number of illustrated examples for the application of the calculus of variations to aerospace problems have been provided in optimal control books, e.g., \cite{Lawden_book, Bryson_book, Longuski_book, Ben_Asher_book, Subchan_book}. 
In this paper, using the calculus of variations, we solve space interception problems with multi-constraints. 

The interception problem of spacecraft is a typical aerospace problem, specifically a space maneuvering problem in which 
an interceptor engaging a target in exoatmosphere 
should 
be maneuvered to approach the target such that
at an impact time,  their position vectors are equal to each other. 
A spacecraft is normally propelled by continuous thrust to realize a space maneuvering. To simplify analysis, a continuous thrust may be approximated by an impulsive thrust when the duration of the thrust is far less than the interval between thrusts or the coasting time of a spacecraft. 
It is also supposed that at a time instant applying an impulsive thrust, the position of a spacecraft is fixed, and its velocity changes by a jump \cite{Curits_book, prussing_2010}. In this sense, we say a velocity impulse instead of an impulsive thrust. Using velocity impulses to approximate continuous thrusts and then obtain an optimal solution is the first step to solve a practical interception problem.   

In 1925,  Hohmann conceived the idea of velocity impulse and first  presented the well-known two-impulse optimal orbit transfer between two circular coplanar space orbits which is now called the Hohmann transfer.  During the 1960s, the subject of impulsive trajectories including interception, transfer, and rendezvous received much attention; see the survey papers \cite{Gobetz1969, Robinson_1967}. 
 By using a variational method, Lawden investigated the optimal control problem of a spacecraft in an inverse square law field 
in which the primer vector theory was specifically developed for impulsive trajectories and a set of necessary conditions 
was presented in \citep{Lawden_book}; see also \cite{Jezewski1975}. 
The primer vector theory was extended to incorporate final time constraints and path constraints for optimal impulsive orbital interception problems by Vinh {\em et al.}~and Taur {\em et al.}~\cite{Vinh1990, Taur1995} respectively.  
The primer vector theory was also applied to multiple-impulse orbit rendezvous, e.g., see \cite{Prussing1986, Sandrik2006, Colasurdo1994}.  For a linear system, Prussing showed the sufficiency of  these necessary conditions for an optimal trajectory \cite{Prussing1995}. 
Recently, a complete design instruction with illustrated numerical examples for the primer vector theory has been presented in \citep{prussing_2010}. 






In our scenario, an interceptor and a target as spacecraft are assumed to be in an inverse square law field without considering other factors, e.g., aerodynamic drags,   
except that the interceptor has two opportunities to adjust its trajectory by using two unknown velocity impulses in order to hit the target and satisfy all constraints simultaneously. 
The initial position and velocity vectors of both interceptors and targets are given. 
Then after given initial states, the interceptor during the period without velocity impulses and the target  both move in free flight. 
In addition to the interception constraint,  the interceptor must satisfy the following multi-constraints. The first kind of constraints is the terminal position constraint of the interceptor at an unknown terminal time, the second is 
the time constraint related to impulse instants at which a velocity impulse is applied and an impact time, 
and the third meaningful constraint is the component-wise magnitude constraint on velocity impulses, which together with time constraints results in two-impulse optimal solutions. With these 
constraints, 
three two-impulse space interception problems are defined.
We then formulate these two-impulse space interception problems as minimum-fuel optimization problems with equality or inequality constraints. 
A dynamic slackness variable method is developed in order to convert the the component-wise inequality constraints on the terminal position of an interceptor into equalities. 
By the calculus of variations,  a set of necessary conditions is established which   finally leads to multi-point boundary value problems with unknown impulse and impact instants and velocity impulses. Therefore, we obtain an indirect optimization method for  two-impulse space interception problems with multi-constraints. 



When ballistic missiles pass through the atmosphere during their flight, they have the same dynamics as spacecraft. 
Hence 
the two-impulse interception problems of free-flight ballistic missiles can be considered as a special kind of two-impulse space interception problems proposed in this paper. 
We then make use of the indirect optimization method developed in this paper to solve the two-impulse interception problems of free-flight ballistic missiles. 
The corresponding boundary value problems have not closed-form analytic solutions, hence we employ the Matlab boundary value problem solvers bvp4c or bvp5c to numerically solve them. 
These two solvers are essentially based on difference method, collocation, and residual control; see \cite{Shampine_book, Kierzenka_thesis, Kierzenka_2008, Kierzenka_2001}, which can deal with multi-point boundary value problems with unknown parameters. 
By using a time change technique, the impulse and impact instants firstly as unknown interior or terminal time instants are introduced as unknown parameters. Then we finally obtain the multi-point boundary value problems with unknown parameters including the components of velocity impulses, which can be solved by the Matlab solvers and are equivalent to the original multi-point boundary value problems directly derived by the calculus of variations. 
A number of conclusions have been established based on the high accuracy of numerical solutions provided by Matlab  solvers.  For example, we find that when time and velocity impulse constraints are imposed, optimal two-impulse  solutions may occur.  

Compared to direct optimization methods, the variational method as an indirect  method provides highly accurate solutions \cite[p.162]{Ben_Asher2014} though the extra initial conditions from costates are needed. This is also verified by numerical examples  in \cite{JOTA_XZX} where for the Hohmann transfer, the numerical solutions provided by the calculus of variations and the Matlab solvers are almost equal to the analytical solutions (only the last three or four digitals of fifteen digitals after decimal place are different).  Hence the calculus of variations is particularly suitable to investigate the properties of numerical solutions.  It is possible that the interception problems could be formulated as static constrained optimization problems since 
the trajectory of a spacecraft in free flight is a conic section. However taking advantages of the calculus of variations and the Matlab  solvers, we formulate them as dynamic constrained  optimization problems.  
In addition, by the dynamics of a spacecraft, we do not need to determine the shape of a conic section (circle, ellipse, parabola, or hyperbola). 
The calculus of variations and the Matlab solvers provide local optimal solutions. 

The organization of this paper is as follows. 
In Section \ref{Sect_1},  three two-impulse space interception problems with equality or inequality constraints are introduced. 
We consider 
the interception problems with equality constraints, i.e., Problems \ref{problem1}-\ref{problem2}  in 
Section \ref{sect_2}. A detailed derivation of necessary conditions for a special case of Problem \ref{problem1} are presented in Appendix A to illustrate how the calculus of variations works, and boundary conditions are given. In order to convert inequality constraints into equalities, two slackness variable methods are presented in Section \ref{sect_3}, and the resulting boundary conditions for Problem \ref{problem3} are obtained.  
The remaining parts of the paper dedicate to solve the two-impulse interception problem of free-flight ballistic missiles.
In Sections \ref{sect_4}-\ref{sect_7}, by using three sets of initial data of the interceptor and target, 
the solution properties of these particular two-impulse space interception problems are characterized and a number of numerical examples are used to illustrate our methods.  
The conclusion is drawn in Section \ref{Sect_8}.  
In Appendix B, three sets of initial data for ballistic missiles are given. 
BCs, parameters, and initial values for Examples \ref{exam1}-\ref{exam2} are shown in Appendix C. 
A time change technique is introduced in In Appendix D. 

  


\section{Problem Statements} \label{Sect_1}

Consider the motion of a spacecraft in the earth-centered inertial frame with the inverse square gravitational field, and the state equation with the force of earth's gravity is 
\begin{align}
\dot{\mathbf r}&=\mathbf v\notag\\
\dot{\mathbf v}&=-\dfrac{\mu}{r^3}\mathbf r\label{eq_Jul_4_01a}
\end{align}
Let an interceptor and a target be such a spacecraft. We only consider the force of earth's gravity and omit other factors. 

\begin{problem} Two-impulse space interception problem \label{problem1}\rm

Given the initial states $\mathbf x_M(t_0)$ and $\mathbf x_T(t_0)$. After the initial time $t_0$, the motion of the target is described by \eqref{eq_Jul_4_01a}. 
Consider two unknown time instants $t_1, t_2$ such that $t_0\leq t_1\leq t_2$.
To guarantee at the impact time $t_h\geq t_2$, the interceptor meets the target, we impose the following equality constraint on the position vectors of the interceptor and the target at an unknown impact time $t_h$ 
\begin{align}\label{eq_Nov5_01f}
\mathbf g_h:=\mathbf{r}_M(t_h)-\mathbf{r}_T(t_h)=0
\end{align}
This is called the interception condition. Suppose that for the interceptor, there are velocity impulses at $t_1$ and $t_2$
\begin{align}
\mathbf{v}_M(t_i^+)=\mathbf{v}_M(t_i^-)+\Delta \mathbf{v}_i, \quad i=1,2
\end{align}
The position vector $\mathbf r(t)$
\begin{align}\label{eq_Nov5_01k}
\mathbf{r}_M(t_i^+)=\mathbf{r}_M(t_i^-), \quad i=1,2
\end{align}
 is continuous at the impulse instants.  
 At the non-impulse instants, the state of the interceptor evolves over time according to the equation \eqref{eq_Jul_4_01a}.  The two-impulse space interception problem as a minimum-fuel optimization problem is to design $\Delta \mathbf{v}_i$ that minimize the cost functional  
\begin{align}
J=\abs{\Delta\mathbf{v}_1}+\abs{\Delta\mathbf{v}_2} \label{eq_Nov14_01}
\end{align}
subject to the constraints \eqref{eq_Nov5_01f}-\eqref{eq_Nov5_01k}. There are three time instants: the impulse instants $t_1,t_2$ and the impact time $t_h$ to be determined.  When there is just only one velocity impulse,  it is called one-impulse space interception problem.
\end{problem}

In control theory, such an optimization problem is called an impulse control problem in which there are state or control jumps.  
Indeed there is not a traditional continuous-time control variable in Problem \ref{problem1}, but instead the unknown velocity impulses can be viewed as  control input parameters defined at discrete time instants.  
It can also be considered as an example of switched or hybrid systems with interior point constraints, e.g., see 
\cite{Bryson_book, Ben_Asher2014}, and references therein.
We next consider a two-impulse space 
interception problem with a terminal position constraint on $\mathbf r_M(t)$. 

\begin{problem}
Two-impulse space interception problem with a terminal position constraint
 \label{problem2}\rm
 
 Consider a similar situation as in Problem \ref{problem1}.  Let  the unknown terminal time  $t_f\geq t_h$. Define a terminal constraint
 \begin{align} \label{eq_Nov5_01l}
\mathbf g_f:=\mathbf r_M(t_f)-\mathbf r_f=0
\end{align}
 on $\mathbf r_M(t_f)$  where $t_f$ is to be determined and the reference position vector $\mathbf r_f=\begin{bmatrix} r_{fx} & r_{fy} & r_{fz}\end{bmatrix}^{\text T}$ is given in exoatmosphere.
The two-impulse space interception problem with a terminal position constraint is to design $\Delta \mathbf{v}_i$ that minimize the cost functional  \eqref{eq_Nov14_01}
subject to the constraints \eqref{eq_Nov5_01f}-\eqref{eq_Nov5_01k} and \eqref{eq_Nov5_01l}. 
\end{problem}

From a practical perspective,  once an interceptor  misses a target for unknown disturbances, the terminal position constraint guarantees 
that the interceptor is in the desired position at a final time. On the other hand,  we will see that by numerical examples, the terminal position constraint can change the impact point of an interception task.
We now propose a two-impulse space interception problem with multi-constraints in terms of inequalities.  

\begin{problem}
Two-impulse space interception problem with multi-constraints
 \label{problem3}\rm
 
 Consider a similar situation as in Problem \ref{problem2}. We have the following multi-constraints:
 \begin{enumerate}
 \item 
 Instead of the equality constraint \eqref{eq_Nov5_01l}, we imposes an inequality constraint on the terminal position $\mathbf r_M(t_f)$
  \begin{align}\label{eq_Nov5_01m}
r_{Mx}(t_f)-r_{fx}-r_{x\max}\leq 0, \quad -\left( r_{Mx}(t_f)-r_{fx}\right)+r_{x\min}\leq 0\notag\\
r_{My}(t_f)-r_{fy}-r_{y\max}\leq 0, \quad -\left( r_{My}(t_f)-r_{fy}\right)+r_{y\min}\leq 0\notag\\
r_{Mz}(t_f)-r_{fz}-r_{z\max}\leq 0, \quad  -\left( r_{Mz}(t_f)-r_{fz}\right)+r_{z\min}\leq 0 
\end{align}
where $\mathbf r_M(t_f)=\begin{bmatrix} r_{Mx} & r_{My} & r_{Mz}\end{bmatrix}^{\text T}$ and the constants 
$r_{x\min}, r_{x\max}, r_{y\min}, r_{y\max},  r_{z\min}, r_{z\max}$ are given. 
\item We also make time constraints on the time instants $t_1, t_2, t_h$
\begin{align}\label{eq_Nov5_01n}
\alpha -t_1\leq 0, \quad t_1-\beta \leq 0, \quad \gamma -(t_2-t_1) \leq 0, \quad \eta-(t_h-t_2)\leq 0
\end{align}
where the time constants $\alpha,\beta, \gamma, \eta\geq 0$ are given. 
\item We impose inequality constraints component-wise on velocity impulses
  \begin{align}\label{eq_Nov5_01o}
\Delta \mathbf v_1=\begin{bmatrix}  v_{1x} &  v_{1y} &  v_{1z} \end{bmatrix}^{\text T}, \quad
\Delta \mathbf v_2=\begin{bmatrix}  v_{2x} &  v_{2y} &  v_{2z} \end{bmatrix}^{\text T}
\notag\\ 
{v_{1x}}\in [p_{1min},  p_{1max}], \quad
{v_{1y}}\in [p_{2min},  p_{2max}], \quad { v_{1z}}\in [p_{3min},  p_{3max}], \notag\\
{v_{2x}}\in [p_{4min},  p_{4max}], \quad
{ v_{2y}}\in [p_{5min},  p_{5max}], \quad { v_{2z}}\in [p_{6min},  p_{6max}]
\end{align}
where all boundary points of the intervals are known.
\end{enumerate}
Then the two-impulse space interception problem with  multi-constraints is to design $\Delta \mathbf{v}_i$ that minimize the cost functional  \eqref{eq_Nov14_01}
subject to the constraints \eqref{eq_Nov5_01f}-\eqref{eq_Nov5_01k} and \eqref{eq_Nov5_01m}-\eqref{eq_Nov5_01o}. 
\end{problem}

Similarly to the terminal position constraint \eqref{eq_Nov5_01l}, the multi-constraints \eqref{eq_Nov5_01m}-\eqref{eq_Nov5_01o} also come from practical requirements.  

\section{Necessary Conditions for Equality Constraints} \label{sect_2}

For simplicity of  presentation, we only derive the first-order necessary condition 
by using the variational method 
for a particular case of Problem \ref{problem1} 
in which the first impulse instant is fixed at the initial time, i.e., $t_1=t_0$. 
Hence the related terms in Problem \ref{problem1} are indexed by $0,1$ instead of $1,2$.

Define the augmented cost functional
\begin{align}\label{eq_Dec6_01}
\tilde J:&=\abs{\Delta\mathbf{v}_0}+\abs{\Delta\mathbf{v}_1}\notag\\&\quad+\mathbf{q}^{\text T}_{v1}\left[\mathbf{v}_M(t_0^+)-\mathbf{v}_M(t_0)-\Delta \mathbf{v}_{0}\right] +\mathbf{q}^{\text T}_{v2}\left[\mathbf{v}_M(t_1^+)-\mathbf{v}_M(t_1^-)-\Delta \mathbf{v}_1\right]
+
\bm\gamma_{h}^\text{T}\mathbf g_{h}(\mathbf r_M(t_h), \mathbf r_T(t_h))\notag\\
&\quad+\int_{t_0^+}^{t_1^-}
\mathbf{p}_{Mr}^{\text{T}}\left (\mathbf{v}_M-\dot{\mathbf{r}}_M\right )+
\mathbf{p}_{Mv}^{\text{T}}\left (-\dfrac{\mu}{r^3_M}\mathbf{r}_M-\dot{\mathbf{v}}_M\right ){\rm d}t
+\int_{t_0^+}^{t_1^-}
\mathbf{p}_{Tr}^{\text{T}}\left (\mathbf{v}_T-\dot{\mathbf{r}}_T\right )+
\mathbf{p}_{Tv}^{\text{T}}\left (-\dfrac{\mu}{r^3_T}\mathbf{r}_T-\dot{\mathbf{v}}_T\right ){\rm d}t
\notag\\
&\quad+\int_{t_1^+}^{t_h}
\mathbf{p}_{Mr}^{\text{T}}\left (\mathbf{v}_M-\dot{\mathbf{r}}_M\right )+
\mathbf{p}_{Mv}^{\text{T}}\left (-\dfrac{\mu}{r^3_M}\mathbf{r}_M-\dot{\mathbf{v}}_M\right ){\rm d}t
+\int_{t_1^+}^{t_h}
\mathbf{p}_{Tr}^{\text{T}}\left (\mathbf{v}_T-\dot{\mathbf{r}}_T\right )+
\mathbf{p}_{Tv}^{\text{T}}\left (-\dfrac{\mu}{r^3_T}\mathbf{r}_T-\dot{\mathbf{v}}_T\right ){\rm d}t
\end{align}
\noindent where the vectors $\mathbf{q}_{v1}, \mathbf{q}_{v2}, \bm\gamma_{h}$ are constant Lagrange multipliers, 
and 
Lagrange multipliers $\mathbf{p}_{Mr}, \mathbf{p}_{Mv}, \mathbf{p}_{Tr}, \mathbf{p}_{Tv}$ are also called costate vectors; in particular, Lawden \cite{Lawden_book} termed $-\mathbf{p}_{Mv}$ the primer vector. 

Introducing Hamiltonian function 
\begin{align}
H_{Mi}(\mathbf{r}_M, \mathbf{v}_M, \mathbf{p}_{Mi}):=\mathbf{p}_{Mri}^{\text{T}}\mathbf{v}_M-\mathbf{p}_{Mvi}^{\text{T}}\dfrac{\mu}{r_M^3}\mathbf{r}_M, \quad\mathbf{p}_{Mi}=\begin{bmatrix}
\mathbf{p}_{Mri} & \mathbf{p}_{Mvi}
\end{bmatrix}^{\text T}, \quad i=1,2\label{eq_7_2_1_Sep19}
\end{align}
Then the integrals  in \eqref{eq_Dec6_01} have simplified expressions, e.g., 
\begin{align}
\int_{t_0^+}^{t_1^-}
\mathbf{p}_{Mr}^{\text{T}}\left (\mathbf{v}_M-\dot{\mathbf{r}}_M\right )+
\mathbf{p}_{Mv}^{\text{T}}\left (-\dfrac{\mu}{r^3_M}\mathbf{r}_M-\dot{\mathbf{v}}_M\right ){\rm d}t:=
\int_{t_0^+}^{t_1^-}
\left(H_{M1}(\mathbf{r}_M, \mathbf{v}_M, \mathbf{p}_{M2})
-\mathbf{p}_{Mr2}^{\text{T}} \dot{\mathbf{r}}_M -
\mathbf{p}_{Mv2}^{\text{T}}\dot{\mathbf{v}}_M\right){\rm d}t\notag
\end{align}

By taking into account of all perturbations, 
we can obtain 
the first variation of the augmented cost functional. After a careful derivation,  from making the coefficients of the variations of all independent variables vanish, 
we derive the boundary conditions related to the costates of the interceptor shown in List \ref{table_Nov7_01}; for a detailed derivation, see Appendix A.  
\vspace*{-.5cm}
\begin{bcs} \label{table_Nov7_01} 
\begin{center} 
\centering{\rm BCs related to the costates for Problem \ref{problem1} \vspace*{0.2cm}}\\
{\blue
\fbox{
\begin{minipage}[H]{7cm} \vspace*{-.3cm}
\black\small
\begin{align}
\begin{array}{ll}
&(1)\begin{cases}
\mathbf p_{Mv1}(t_0^+)+\dfrac{\Delta \mathbf{v}_0}{\abs{\Delta \mathbf{v}_0}}=0, \quad 
\mathbf p_{Mv1}(t_1^-)+
\dfrac{\Delta \mathbf{v}_1}{\abs{\Delta \mathbf{v}_1}} =0\notag\\
\mathbf p_{Mv1}(t_1^-)-\mathbf p_{Mv2}(t_1^+)=0, \quad
\mathbf p_{Mv2}(t_h)=0 \notag\\
\mathbf p_{Mr1}(t_1^{-})-\mathbf p_{Mr2}(t_1^{+})=0\notag
\end{cases}
\vspace*{0.1cm}\\
&(2)~ 
-\mathbf p_{Mr1}^{\text T}(t_1^{-})\Delta \mathbf v(t_1)=0, \quad
\mathbf p_{Mr2}^{\text T}(t_h)\left(\mathbf v_M(t_h)-
\mathbf v_T(t_h)\right)=0\notag\\
\end{array}
\end{align}
\end{minipage}}}
\end{center}
\end{bcs}

The second BCs in List \ref{table_Nov7_01} are from the Hamiltonian conditions. 
The similar arguments can be applied to obtain 
the boundary conditions related to the costates of the interceptor for Problem \ref{problem2} which is shown in List \ref{table_Nov7_02}.
\begin{bcs} \label{table_Nov7_02}
\begin{center} 
\centering{\rm BCs related to the costates for Problem \ref{problem2} \vspace*{0.2cm}}\\
{\blue
\fbox{
\begin{minipage}[H]{7.15cm} \vspace*{-.3cm}
\black\small
\begin{align}
\begin{array}{ll}
&(1)\begin{cases}
\mathbf p_{Mv1}(t_1^-)+\dfrac{\Delta \mathbf{v}_1}{\abs{\Delta \mathbf{v}_1}}=0, \quad \mathbf p_{Mv2}(t_2^-)+
\dfrac{\Delta \mathbf{v}_2}{\abs{\Delta \mathbf{v}_2}} =0\notag\\
\mathbf p_{Mv1}(t_1^-)-\mathbf p_{Mv2}(t_1^+)=0, \quad
\mathbf p_{Mv2}(t_2^{-})-\mathbf p_{Mv3}(t_2^{+})=0 
\notag\\
\mathbf p_{Mr1}(t_1^{-})-\mathbf p_{Mr2}(t_1^{+})=0, \quad 
\mathbf p_{Mr2}^{\text T}
(t_2^{-})-\mathbf p_{Mr3}^{\text T}(t_2^{+})=0\notag\\
\mathbf p_{Mv3}(t_h^{-})-\mathbf p_{Mv4}(t_h^{+})=0, \quad \mathbf p_{Mv4}(t_f)=0
\notag
\end{cases}
\vspace*{0.1cm}\\
&(2)\begin{cases}
-\mathbf p_{Mr1}^{\text T}(t_1^{-})\Delta \mathbf v(t_1)=0, \quad
-\mathbf p_{Mr2}^{\text T}(t_2^{-})\Delta \mathbf v(t_2)=0\notag\\
\left(\mathbf p^{\text T}_{Mr3}(t_h^{-})
-
\mathbf p^{\text T}_{Mr4}(t_h^{+})\right)
\left(\mathbf v_M(t_h)-\mathbf v_T(t_h)\right)=0
\notag\\
 \mathbf p_{Mr4}^{\text T}(t_f)\mathbf v_{M}(t_f)
 =0\notag
\end{cases}
\vspace*{0.1cm}\\
\end{array}
\end{align}
\end{minipage}}}
\end{center}
\end{bcs}


\section{Necessary Conditions for Inequality Constraints}  \label{sect_3}

We are now in a position to give the first-order necessary condition for Problem \ref{problem3} with multi-constraints in terms of inequalities.  The inequality constraints must be first converted into equality constraints in order to use Lagrange multiplier method.  We use two kinds of slackness variable methods, one of which is rooted in the well-known Kuhn-Tucker Theorem, and the other is a dynamic slackness variable method developed in this paper. 

\subsection*{Inequality constraints on time instants and velocity impulses} \label{Sect_4_1}

For the inequality constraints \eqref{eq_Nov5_01n}-\eqref{eq_Nov5_01o}, we only consider the first inequality  $\alpha-t_1\leq 0$ in the time constraints \eqref{eq_Nov5_01n} for a simple presentation. We use Problem \ref{problem1} to illustrate how to convert an inequality constraint into an equality constraint in light of nonlinear programming. 

\medskip
\noindent {\emph{A static slackness variable method}}

\medskip

We now introduce the slackness variable $\epsilon$ which ensures that $t_1\geq \alpha$
\begin{align}
t_1-\alpha-\epsilon^2=0 \Longrightarrow t_1-\alpha =\epsilon^2\geq 0\label{eq_Oct3_01}
\end{align}
Then we obtain an equality constraint which can be dealt with by Lagrange multiplier method; e.g., see \cite[Section 4.3]{Hull_book} for slackness variable. Define a new augmented cost function
\begin{align}
\hat J=\tilde J+\lambda(t_1-\alpha -\epsilon^2)\label{eq_Dec27_01}
\end{align}
where $\tilde J$ is the original augmented cost functional defined by \eqref{eq_Dec6_01}, and $\lambda$ is an unknown constant Lagrange multiplier. 
We then calculate the first variation of $\hat J$
\[
\delta\hat J=\delta\tilde J+\lambda\,\delta t_1-2\lambda \epsilon \,\delta\epsilon
\]
In order to make $\delta \epsilon$ vanish such that $\delta\hat J=0$, we let 
$\lambda\epsilon=0$, and meanwhile group $\lambda\delta t_1$
into the related terms of $\delta t_1$ in $\delta\tilde J$, for example, in Problem \ref{problem1}, we have 
\begin{align}
-\mathbf p_{Mr1}^{\text T}(t_1^{-})\Delta \mathbf v(t_1)+\lambda=0 \label{eq_Oct3_02}
\end{align} 
 Eliminating the slack variable yields 
\begin{align}
\lambda\epsilon=0\Longrightarrow \lambda\epsilon^2=0 \stackrel{\eqref{eq_Oct3_01}}{\Longrightarrow} \lambda (t_1-\alpha)=0 \label{eq_Oct3_03}
\end{align}
We now obtain the necessary condition for the inequality constraint $\alpha-t_1\leq 0$ in terms of the interior point boundary conditions \eqref{eq_Oct3_02}-\eqref{eq_Oct3_03} in which the Lagrange multiplier $\lambda$ is an unknown constant to be determined.  A similar argument can be used to the case of component-wise inequality constraints on velocity impulses.  We find that our method is actually  a version of the Kuhn-Tucker Theorem in which \eqref{eq_Oct3_03} is called complementary slackness condition. However the condition \eqref{eq_Oct3_02} is new and attributes to the variational method to the two-impulse interception problem \ref{problem1}. 

\subsection*{Inequality constraints on the position vector of the interceptor}

Considering the terminal position constraints \eqref{eq_Nov5_01m} on $\mathbf r(t_f)$, instead of the static slackness variable method, 
we develop a dynamic slackness variable method.  We next only consider the constraint on the $x$-direction component $r_{Mx}(t_f)$ of the position vector $\mathbf r_M(t_f)$ to just illustrate how the dynamic slackness variable method works. 

\medskip
\noindent {\emph{A dynamic slackness variable method}}


\medskip
We introduce the dynamic slackness variable $\epsilon(t)$ by the left-hand side of \eqref{eq_Oct19_03}
\begin{equation}
\dot \epsilon=
\begin{cases}
0, & t\in [t_0, t_1)\cup [t_1, t_2)\\
w, & t\in [t_2, t_h)\\
w, & t\in [t_h, t_f]
\end{cases}, \qquad
\dot p_\epsilon=
\begin{cases}
0, & t\in [t_0, t_1)\cup [t_1, t_2)\\
-2k_3\epsilon, & t\in [t_2, t_h)\\
-2k_4\epsilon, & t\in [t_h, t_f]
\end{cases}
\label{eq_Oct19_03}
\end{equation}
where $w$ is an unknown variable to be optimized. 
Using $\epsilon(t)$, we convert the first two inequalities in \eqref{eq_Nov5_01m} into the equality constraints as follows
\begin{align}
r_{Mx}(t_f)-r_{fx}-r_{x\max}+\epsilon^2(t_f)=0, \quad -\left( r_{Mx}(t_f)-r_{fx}\right)+r_{x\min}+\epsilon^2(t_h^-)=0\notag
\end{align}
By adding two integrals of quadratic forms
in the slackness variable $\epsilon$ and its control $w$ and also Lagrange multiplier terms into the original cost functional $\tilde J$, we have the augmented cost functional $\hat J$
\begin{align}
\hat J&=\tilde J+\int_{t_1^+}^{t_h^-}\left(w^2+k_3\epsilon^2\right)d t+\int_{t_h^+}^{t_f}\left(w^2+k_4\epsilon^2\right)d t\notag\\
&\quad 
+\eta_1\left(r_{Mx}(tf)-r_{fx}-r_{x\max}+\epsilon^2(t_f)\right) +\eta_2\left(-\left( r_{Mx}(tf)-r_{fx}\right)+r_{x\min}+\epsilon^2(t_h^-)
\right)\label{eq_Dec15_01}
\end{align}
where $k_3, k_4$ are weighting coefficients.  By the calculus of variations, we now solve this dynamic optimization problem in which the dynamical equation consists of \eqref{eq_Jul_4_01a} both for the interceptor and the target and also combines with the slackness variable equation \eqref{eq_Oct19_03}. 
We have the optimal $w=-0.5p_\epsilon$ and the costate equation given by the right-hand side of \eqref{eq_Oct19_03}. 
In the first order variation of $\hat J$,  one can easily obtain the terms related to the variation of $\epsilon$ as follows
\begin{align}
p_{1\epsilon}(t_0)\delta \epsilon(t_0)
-p_{1\epsilon} (t_1^{-*})\delta \epsilon(t_1^{-*})
+p_{2\epsilon} (t_1^{+*})\delta \epsilon(t_1^{+*}),\quad
-p_{2\epsilon} (t_2^-)\delta \epsilon(t_2^-)
+p_{3\epsilon} (t_2^{+*})\delta \epsilon(t_2^{+*}),
\notag\\
-p_{3\epsilon} (t_h^{-*})\delta \epsilon(t_h^{-*})
+p_{4\epsilon} (t_h^{+*})\delta \epsilon(t_h^{+*})-p_{4\epsilon} (t_f^*)\delta \epsilon(t_f^*)\label{eq_May6_01}
\end{align}
where we  denote the costate $p_{\epsilon}$ in the $i$-{\em th} time interval as  $p_{i\epsilon}, i=1,\dots,4$.
Using a similar equality as \eqref{eq_Oct19_eq4} in Appendix A, 
eliminating the variations $\delta \epsilon$ and 
regrouping the terms in $d\epsilon$, 
 we have the boundary conditions for $\epsilon$ and $p_\epsilon$ 
\begin{align}
p_{1\epsilon}(t_0)=0, \quad p_{1\epsilon}^{\text T}(t_1^{-*})-p_{2\epsilon}^{\text T}(t_1^{+*})=0, \quad
p_{2\epsilon}^{\text T}(t_2^{-*})-p_{3\epsilon}^{\text T}(t_2^{+*})=0\notag\\
-p_{3\epsilon}^{\text T}(t_h^{-*})
+p_{4\epsilon}^{\text T}(t_h^{+*})+2\eta_2\epsilon(t_h^{-*})=0, \quad
-p_{4\epsilon}^{\text T}(t_f^*)+2\eta_1\epsilon(t_f^*)=0\notag\\
\epsilon(t_1^{-*})-\epsilon(t_1^{+*})=0, \quad
\epsilon(t_2^{-*})-\epsilon(t_2^{+*})=0, \quad
\epsilon(t_h^{-*})-\epsilon(t_h^{+*})=0\notag
\end{align}
where the first five terms come from making the coefficients of $d\epsilon$  vanish, and the last three terms follow from the continuous of $\epsilon$.  
Notice that the relationship between the Hamiltonian function $\hat H$ in the new augmented system and 
the original $H$ are as follows
\begin{align}
\hat H_i=H_i, \quad i=1,2; \quad
\hat H_i=H_i+w^2+p_{i}w+k_{i}\epsilon^2=H_i-0.25p_{i}^2+k_{i}\epsilon^2, \quad i=3,4\label{eq_Dec27_02}
\end{align}
Applying the same arguments to the $y,z$-direction components of $\mathbf r_M(t_f)$, we finally obtain a complete boundary conditions related to the costates and inequality constraints for multi-constraints in List \ref{table_Oct21_01} 
in which
\begin{align}
\bm \epsilon=\begin{bmatrix} \epsilon_x & \epsilon_y & \epsilon_z\end{bmatrix}^{\text T}, \quad \bm \epsilon^2=\begin{bmatrix} \epsilon_x^2 & \epsilon_y^2 & \epsilon_z^2\end{bmatrix}^{\text T}, \quad \mathbf p_{\epsilon}=\begin{bmatrix} p_x & p_y & p_z\end{bmatrix}^\text{T}\notag\\
\mathbf r_{\min}=\begin{bmatrix} r_{x\min} & r_{y\min} & r_{z\min}\end{bmatrix}^\text{T}, \quad
\mathbf r_{\max}=\begin{bmatrix} r_{x\max} & r_{y\max} & r_{z\max}\end{bmatrix}^\text{T}, \quad
\Delta \hat H_1=\hat H_3- H_3, \quad \Delta \hat H_2=\hat H_4- H_4\label{eq_Dec12_01}
\end{align}
\newpage
\begin{bcs} \label{table_Oct21_01}
\begin{center} 
\centering{\rm BCs related to the costates and inequality constraints\vspace*{0.2cm}}
{\blue
\fbox{\black
\begin{minipage}[H]{1cm} \vspace*{-.3cm}
\small
\begin{align}
\begin{array}{ll}
\vspace*{0.2cm}
&(1)\begin{cases}
\mathbf p_{Mv1}(t_1^-)+\dfrac{\Delta \mathbf{v}_1}{\abs{\Delta \mathbf{v}_1}}+
\begin{bmatrix} \mu_1-\mu_2 & \mu_3-\mu_4 &  \mu_5-\mu_6\end{bmatrix}^{\text T}
=0\notag\\
\mathbf p_{Mv2}(t_2^-)+
\dfrac{\Delta \mathbf{v}_2}{\abs{\Delta \mathbf{v}_2}}  
+
\begin{bmatrix} \mu_7-\mu_8 & \mu_9-\mu_{10} &  \mu_{11}-\mu_{12}\end{bmatrix}^{\text T}
=0\notag\\
\mathbf p_{Mv1}(t_1^-)-\mathbf p_{Mv2}(t_1^+)=0, \quad
\mathbf p_{Mv2}(t_2^{-})-\mathbf p_{Mv3}(t_2^{+})=0 
\notag\\
\mathbf p_{Mv3}(t_h^{-})-\mathbf p_{Mv4}(t_h^{+})=0, \quad
 \mathbf p_{Mv4}(t_f)=0
\notag\\
\mathbf p_{Mr1}(t_1^{-})-\mathbf p_{Mr2}(t_1^{+})=0, \quad
\mathbf p_{Mr2}^{\text T}
(t_2^{-})-\mathbf p_{Mr3}^{\text T}(t_2^{+})=0
\notag
\end{cases}\\
\vspace*{0.2cm}
&(2)\begin{cases}
-\mathbf p_{Mr1}^{\text T}(t_1^{-})\Delta \mathbf v(t_1)-\lambda_1+\lambda_2+\lambda_3=0\notag\\
-\mathbf p_{Mr2}^{\text T}(t_2^{-})\Delta \mathbf v(t_2)-\lambda_3+\lambda_4=0 \notag\\
\left(\mathbf p^{\text T}_{Mr3}(t_h^{-})
-
\mathbf p^{\text T}_{Mr4}(t_h^{+})\right)
\left(\mathbf v_M(t_h)-\mathbf v_T(t_h)\right)-\lambda_4+\Delta\hat H_1=0
\notag\\
 \mathbf p_{Mr4}^{\text T}(t_f)\mathbf v_{M}(t_f)+\Delta\hat H_2
 =0\notag
\end{cases}\\\vspace*{0.2cm}
&(3)\begin{cases}{{\lambda_1(\alpha-t_1)=0, \quad \lambda_2(t_1-\beta)=0}}\notag\\
\lambda_3(\gamma-(t_2-t_1))=0, \quad \lambda_4(\eta-(t_h-t_2))=0\notag
\end{cases}\\\vspace*{0.2cm}
&(4)\begin{cases}\vspace*{0.1cm}
\diag(\mu_1, \mu_3,  \mu_5)
\left(\Delta \mathbf{v}_1-\mathbf p_{1\max}\right)=0
\\\vspace*{0.1cm}
\diag(\mu_2, \mu_4, \mu_6)
\left(\mathbf p_{1\min}
-\Delta \mathbf v_1\right)=0\\\vspace*{0.1cm}
\diag(\mu_7,\mu_9,\mu_{11})
\left(\Delta \mathbf{v}_2-\mathbf p_{2\max}\right)=0
\\\vspace*{0.1cm}
\diag(\mu_8, \mu_{10}, \mu_{12})
\left(\mathbf p_{2\min}
-\Delta \mathbf v_2\right)=0
\end{cases}\\\vspace*{0.2cm}
&(5)\begin{cases}
\bm \epsilon(t_1^+)-\bm \epsilon(t_1^-)=0, \quad 
\bm \epsilon(t_2^+)-\bm \epsilon(t_2^-)=0, \quad
\bm \epsilon(t_h^+)-\bm \epsilon(t_h^-)=0\\
\mathbf p_\epsilon(t_1^+)-\mathbf p_\epsilon(t_1^-)=0, \quad 
\mathbf p_\epsilon(t_2^+)-\mathbf p_\epsilon(t_2^-)=0 \\
\end{cases}\\\vspace*{0.2cm}
&(6)\begin{cases}
-\mathbf p_\epsilon(t_h^-)+\mathbf p_\epsilon(t_h^+)+2
\diag(\eta_2,\eta_4,\eta_6)\bm \epsilon(t_h^-)=0\\
-\mathbf p_\epsilon(t_f)+2
\diag(\eta_1,\eta_3,\eta_5)\bm \epsilon(t_f)=0, \quad \mathbf p_{\epsilon}(t_0)=0\\
\end{cases}
\\
&(7)\begin{cases}
\left(\mathbf r_{M}(t_f)-\mathbf r_f\right)-\mathbf r_{\max}+\bm \epsilon^2(t_f)=0\\
-\left(\mathbf r_{M}(t_f)-\mathbf r_f\right)+\mathbf r_{\min}+\bm \epsilon^2(t_h^-)=0\\
\end{cases}
\\
\vspace*{0.4cm}
&(8)~~
\begin{bmatrix}
\eta_1-\eta_2 &
\eta_3-\eta_4 &
\eta_5-\eta_6 
\end{bmatrix}^{\text T}
-\mathbf r(t_f)=0
\end{array}
\notag 
\end{align}
\end{minipage}}}
\end{center}
\end{bcs}

To obtain the BCs in List \ref{table_Oct21_01}, for the time constraints and velocity impulse constraints, we use the static slackness variable method; conversely, for the terminal position inequality constraints, the dynamic slackness variable method works. 
The first group describes the optimal conditions of the velocity impulses with component-wise constraints  and the continuous conditions of the costates.
The first two boundary conditions
in the second group is  due to Hamiltonian functions of the original intercept problem \ref{problem1} and the time inequality constraints \eqref{eq_Nov5_01n}, and $\Delta \hat H_1, \Delta \hat H_2$ defined by \eqref{eq_Dec12_01} in the last two boundary conditions are associated with the slackness variables associated with time instants $t_h^-$ and $t_f$.
The boundary conditions in the third and fourth group are complementary slackness conditions similar to \eqref{eq_Oct3_03}, associated with the time inequality constraints \eqref{eq_Nov5_01n} and component-wise magnitude constraints on velocity impulses. 
The remaining groups come from the terminal position component-wise constraints and also the continuous conditions of the dynamic slackness variable and its costate.

\section{One-impulse Space Interception Problem} \label{sect_4}

Starting from this section, by using all derived BCs, we now investigate the properties of two-impulse space interception problems of free-flight ballistic missiles.  

We first consider the one-impulse space interception problem in which there is one velocity impulse to realize an interception task. 
When we deal with two-impulse interception problems by using Matlab boundary value problem solvers, for our initial data generated by ballistic missiles, we find that if the first impulse instant is fixed at the initial time and the second one is free, then its optimal solution degenerates to the one-impulse case in which the optimal impulse instant is equal to the initial time, which means that a one-impulse solution is optimal.  Therefore, it is necessary to study the one-impulse space interception problem first.


Without loss of generality, we assume $t_0=0$. 
Let the impulse instant $t_1$ be fixed. We have $\delta t_1=0$ in the first variation of the cost $\tilde J$ \eqref{eq_Dec6_01}. Hence the first boundary condition in the second group of List \ref{table_Nov7_01} related to $t_1$ vanishes. Then we have BCs for the one-impulse space interception problem with fixed impulse instants as follows
\begin{align}
&(1)\begin{cases}
\mathbf p_{Mv1}(t_1^-)+\dfrac{\Delta \mathbf{v}_1}{\abs{\Delta \mathbf{v}_1}}=0\notag\\
\mathbf p_{Mv1}(t_1^-)-\mathbf p_{Mv2}(t_1^+)=0, \quad
\mathbf p_{Mv2}(t_h)=0 \notag\\
\mathbf p_{Mr1}(t_1^{-})-\mathbf p_{Mr2}(t_1^{+})=0
\end{cases}\\
&(2)~
\mathbf p^{\text T}_{Mr2}(t_h)\left(\mathbf v_M(t_h)-\mathbf v_T(t_h)\right)=0\notag
\end{align}
Notice that the first BC in the first group of List \ref{table_Nov7_01}  has been removed since we consider the one-impulse case.  

\begin{figure}
 \begin{minipage}{8cm}
 \centering
 \includegraphics[width=\hsize]{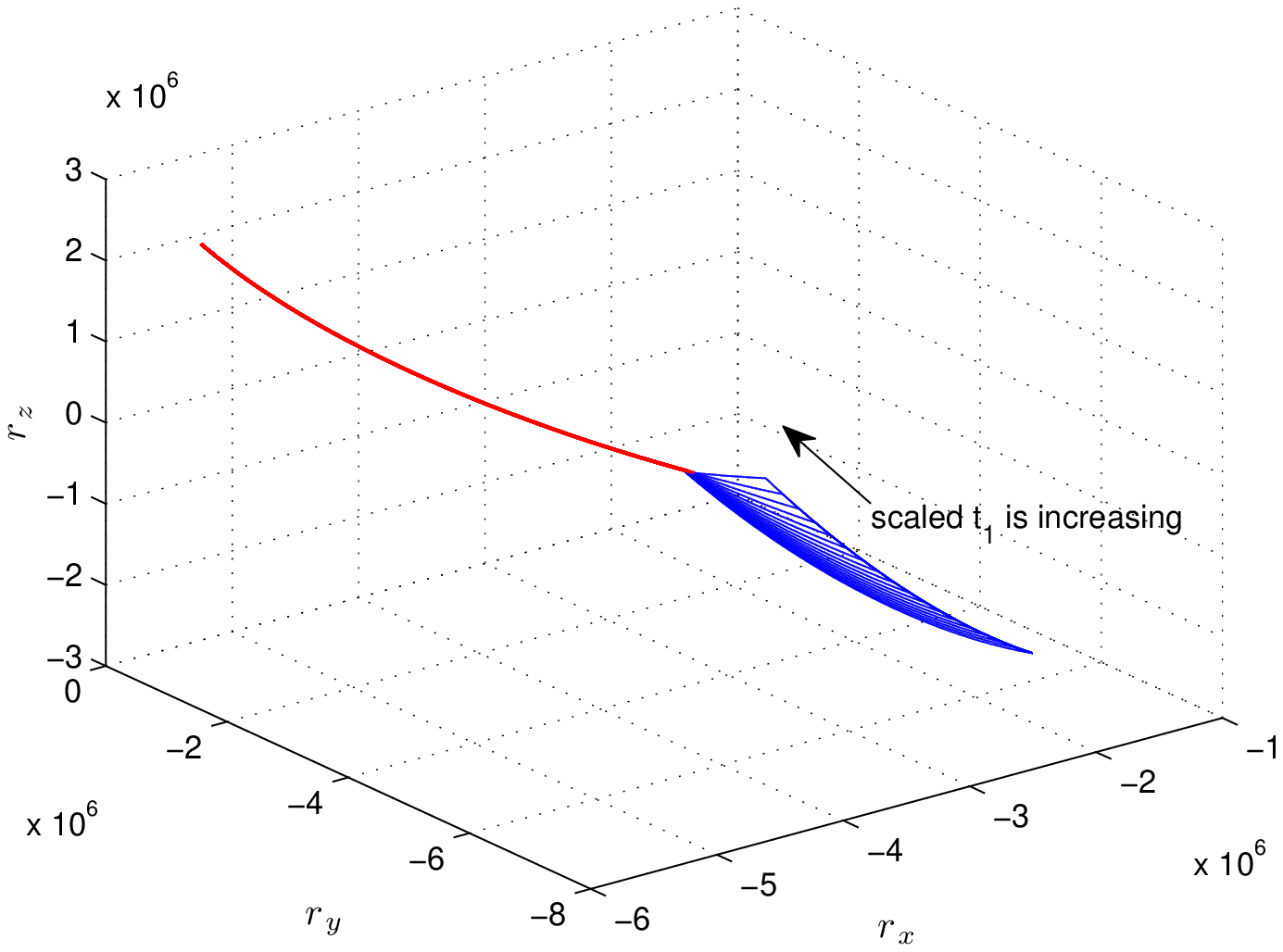}
\end{minipage}
\qquad 
 \begin{minipage}{8cm}
 \centering
\includegraphics[width=\hsize]{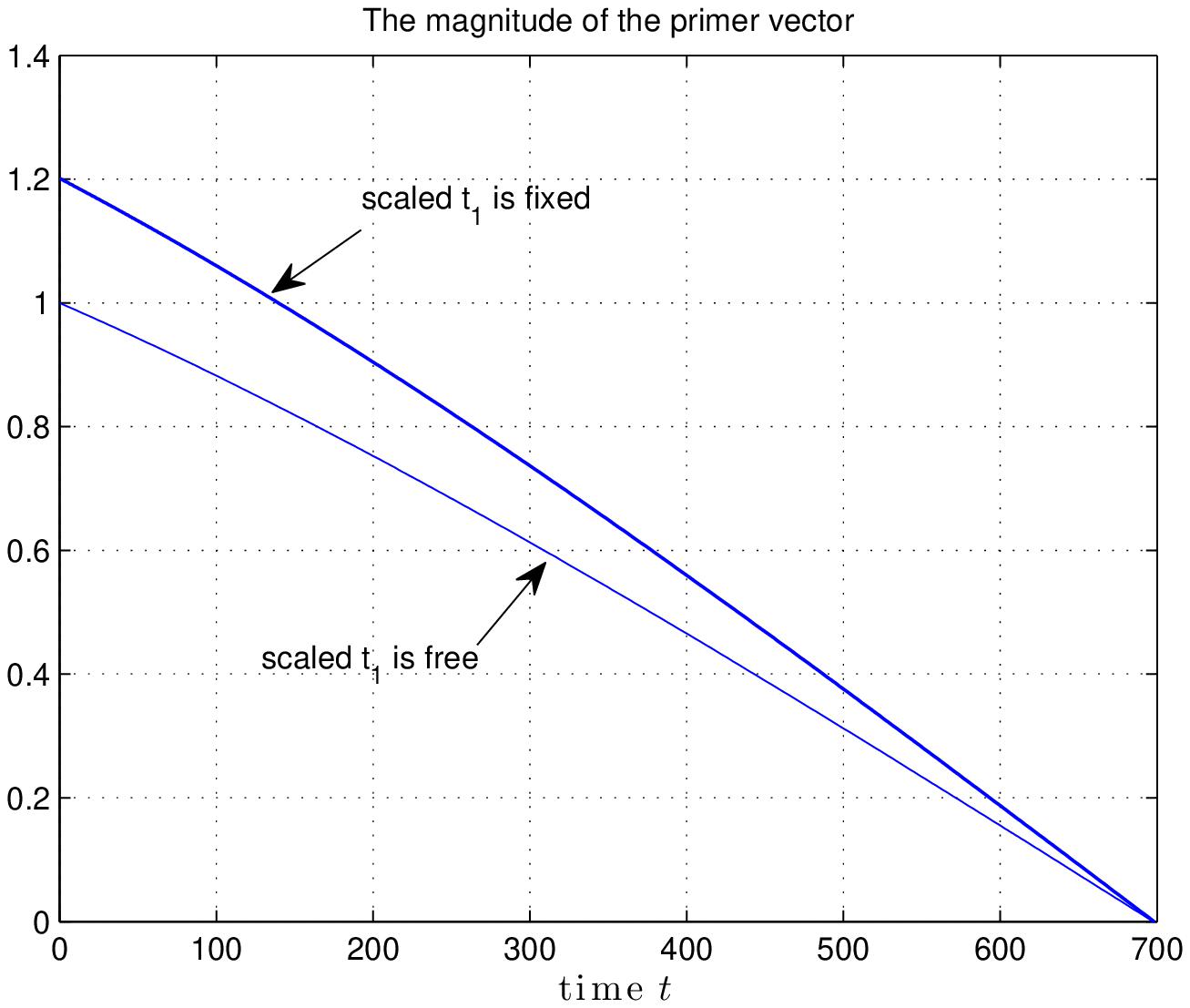}
\end{minipage}
\caption{
Interception trajectories vs.~scaled $t_1$ and the primer vector magnitude for scaled $t_1=0.2$}\label{fig_exam3_4_911_1}
\end{figure}

\begin{figure}[h!]
 \begin{minipage}{8cm}
 \centering
 \includegraphics[width=\hsize]{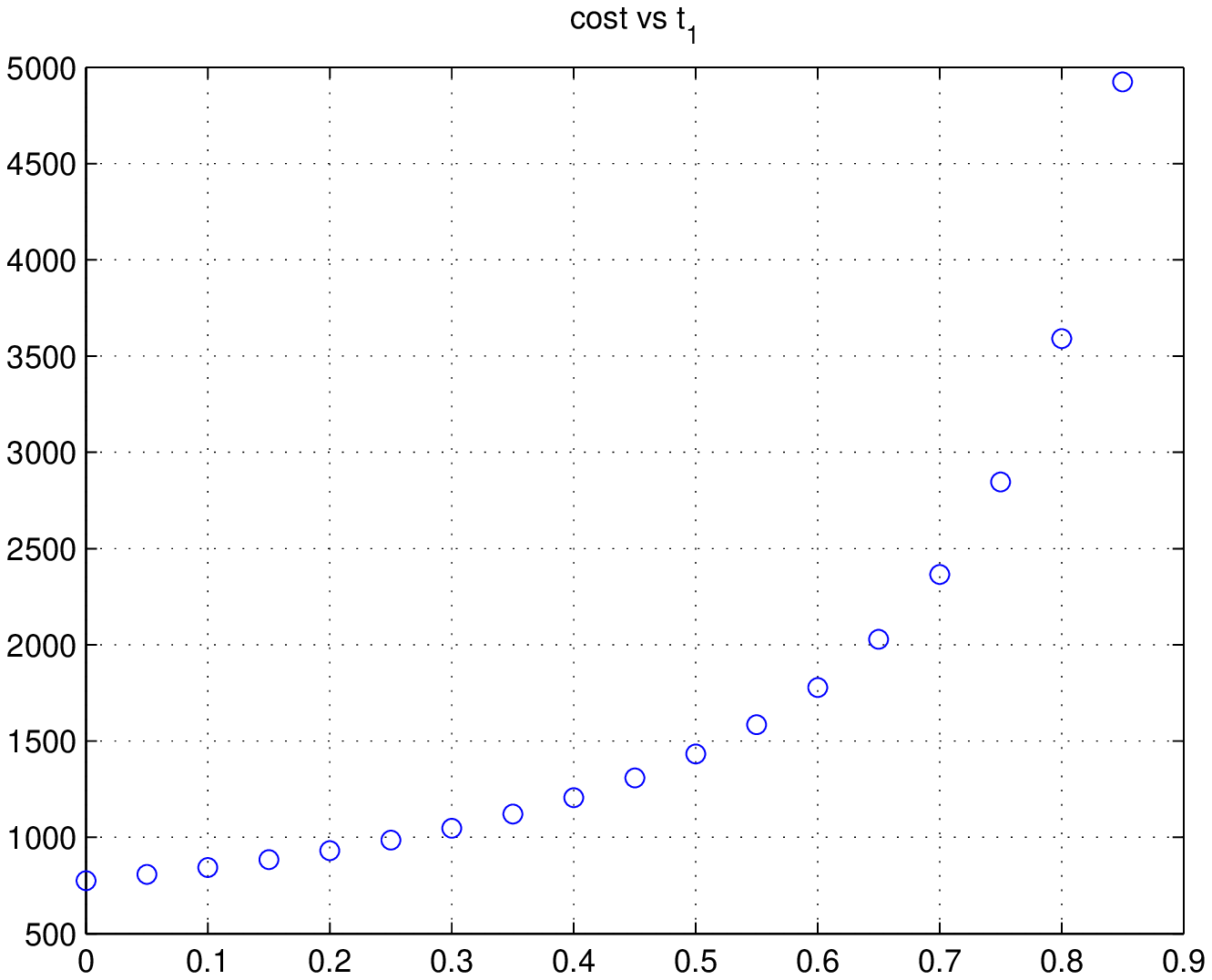}
\caption{Cost vs. scaled impulse instant} \label{fig_exam3_3_a_1}
\end{minipage}
\qquad 
 \begin{minipage}{8cm}
 \centering
 \includegraphics[width=\hsize]{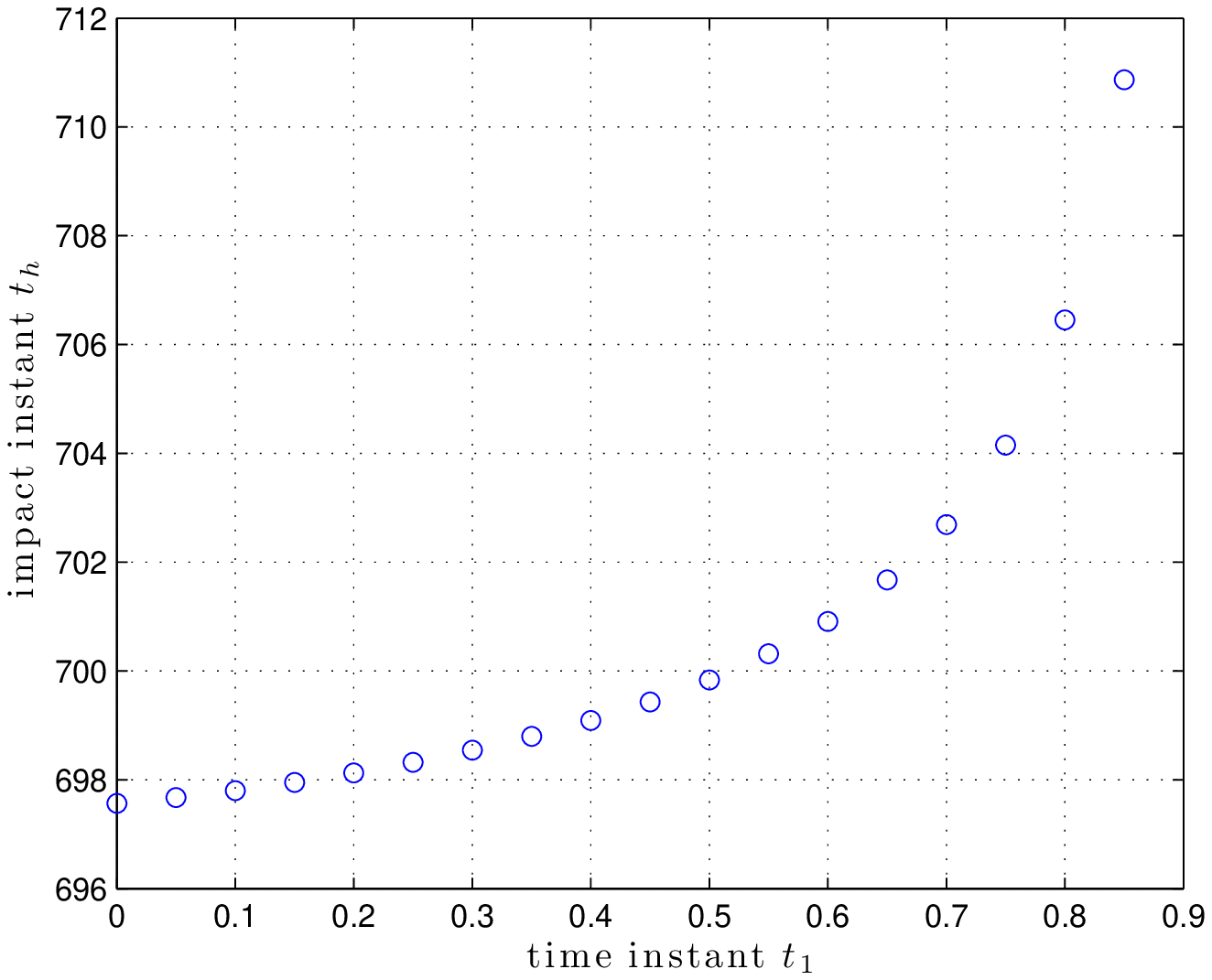}
\caption{Impact instant vs. scaled impulse instant}\label{fig_exam3_4_b}
\end{minipage}
\end{figure}

Three sets of initial data numerically generated by the elliptical orbits of ballistic missiles as interceptors and targets have been provided in Appendix B. 
In order to use the Matlab solvers, the impulse instant $t_1$ is normalized by the impact instant $t_h$ 
after a time change is introduced; see Appendix C for the definition of scaled time instants. We call such a normalized instant a scaled instant.  Here 
the impulse instant $t_1$ equals $\text{scaled~} t_1\times t_h$.
Local optimal numerical solutions for each scaled $t_1=[0:0.05:0.85]$ have been obtained for the initial data I.  
The left figure of Fig.~\ref{fig_exam3_4_911_1}
shows interception trajectories for each scaled $t_1=[0:0.05:0.85]$
in which the upward direction of trajectories is the direction of increasing time.
The primer vector magnitude for the fixed scaled $t_1=0.2$ is shown in the right figure (the top one). One can see that its magnitude at the impulse instant $t_1$ is equal to one and does not achieve its maximum 
which does not satisfy the necessary condition of optimal transfers given by Lawden \cite{Lawden_book} (see also Table 2.1 Impulsive necessary conditions in \citep{prussing_2010}),
hence this impulse instant $t_1$ is not optimal.  
The magnitude of velocity impulse  and the impact instant $t_h$ vs.~scaled impulse instant $t_1$ are shown in Fig.~\ref{fig_exam3_3_a_1}-\ref{fig_exam3_4_b} respectively. 
Obviously, with the increasing of scaled $t_1$, the magnitude of velocity impulse  strictly increases, which is consistent with our intuition.


Table \ref{table03} gives the specific data of fixed impulse instants, impact instants, and the amplitudes of the velocity impulse obtained by the solver bvp5c in which we set Tolerance=1e-9. It is found that the larger the impulse instant $t_1$, the more time the computation takes. We now conclude the above numerical results by Claim \ref{claim1}.

\begin{conclusion}\label{claim1}
\rm For the initial data I, 
the optimal cost of the one-impulse space interception problem is monotonically increasing with the impulse instant $t_1$.
\end{conclusion}

\begin{center}
\begin{threeparttable}\small 
\caption{Fixed impulse instants, impact instants, and velocity impulse magnitudes}\label{table03}
\begin{tabular}{ccccccccc}
\hline
$t_1$ & $t_h$ & $\abs{\Delta\mathbf v_1}$ &  $t_1$ & $t_h$ & $\abs{\Delta\mathbf v_1}$\\
\hline
0 & 697.5637 & 774.9142 & 60 & 697.7622 & 832.3724\\
5  & 697.5779 & 779.2491 & 70 & 697.8013 & 843.2137\\
10 &697.5925 &783.6609 & 80 & 697.8422 & 854.4603\\
20 & 697.6230&792.7212&90 & 697.8847 & 866.1314\\
30 & 697.6553 & 802.1083 & 100 & 697.9291 & 878.2470\\
40 & 697.6892& 811.8358 & 200 & 698.4787 & 1.0295e+3\\
50 & 697.7249&821.9186  & 600 & 710.2164 & 4.7368 e+3\\ 
\hline
\end{tabular}
\end{threeparttable}
\end{center}
\medskip

It follows from Table \ref{table03} as well as Claim \ref{claim1} that  the cost of the one-impulse space interception problem attains the minimum at the initial time $t_0=0$. This can also be seen from solving the related boundary value problem with the free instant $t_1$ by imposing an inequality constraints $-t_1\leq 0$ for the one-impulse case.  Without this inequality constraint, the Matlab solver gives a negative impulse instant. It is reasonable since in the derivation of BCs by the variational method, we do not impose the time constraint $t_1\geq 0$. 
Notice that the inequality $t_1\geq 0$ is essentially not a constraint since the time is naturally non-negative. Hence for this case we also say that the instant $t_1$ is free. 
List \ref{table_Oct5_01} shows the BCs by using the static slackness variable method; see \eqref{eq_Oct3_02}-\eqref{eq_Oct3_03}. 

\begin{conclusion}\label{claim2}
\rm For the initial data I, 
the one-impulse space interception problem has the minimum at the initial time $t_0$. 
\end{conclusion}

\begin{bcs} \label{table_Oct5_01}
\begin{center} 
\centering{\small BCs related to the costaes for free $t_1$\vspace*{0.2cm}}

{\blue
\fbox{
\begin{minipage}[H]{5cm} \vspace*{-.3cm}
\black\small
\begin{align}
\begin{array}{ll}
&(1)\begin{cases}
 \mathbf p_{Mv1}(t_1^-)+
\dfrac{\Delta \mathbf{v}_1}{\abs{\Delta \mathbf{v}_1}} =0\\
\mathbf p_{Mv1}(t_1^-)-\mathbf p_{Mv2}(t_1^+)=0, \quad
\mathbf p_{Mv2}(t_h )=0\\
\mathbf p_{Mr1}(t_1^{-})-\mathbf p_{Mr2}(t_1^{+})=0
\end{cases}
\vspace*{0.1cm}
\\
&(2)\begin{cases}
-\mathbf p_{Mr1}^{\text T}(t_1^{-})\Delta \mathbf v(t_1){-\lambda}=0\\
 \mathbf{p}_{Mr2}(t_h) \left(\mathbf{v}_{M}(t_h)-\mathbf{v}_{T}(t_h)\right)=0
\end{cases}\vspace*{0.1cm}\\
&(3)~{{\lambda (0-t_1)=0}}\\
\end{array}
\notag
\end{align}
\end{minipage}}}
\end{center}
\end{bcs}

\bigskip
By the BCs in List \ref{table_Oct5_01}, the solver bvp5c with Tolerance 2.22045e-14 gives the following solution message
{\small
\begin{verbatim}
The maximum error is 2.068e-14. Elapsed time is 132.710440 seconds.
dv1 = [-376.7263; 338.2633; -586.6407], abs(dv1) = 774.9142
t_impulse = 0, t_impact = 697.5637, lambda = 0.8594
\end{verbatim}
}
Then one can see that the local optimal impulse instant exactly equals zero, which together with Claim \ref{claim1} numerically verifies Claim \ref{claim2}.  
The primer vector magnitude shown in the right figure (the bottom one)  achieves its maximum one at $t=0$, which is consistent with the necessary condition of the primer vector theory; see item 2 of Table 2.1 in \cite{prussing_2010}.
\begin{rem} \rm
The numerical experiments for each fixed (scaled) impulse instant $t_1$ in this section are carried out under the same initial circumstance including the same initial values of state and parameters, and Claims \ref{claim1}-\ref{claim2} are based on the first-order necessary condition for optimality derived from the calculus of variations and the Matlab solvers both of which can guarantee only the local optimality of solutions.
In this sense, by Claims \ref{claim1}-\ref{claim2}, we report that the local optimal cost for the one-impulse interception problem of ballistic missiles has a monotonic property  with respect to the varying impulse instants under the same initial numerical circumstance, and takes a minimum at the initial instant. 
In the sequel, a local optimal solution is quite often simply called an optimal solution. 
It is noted that Claims \ref{claim1}-\ref{claim2} also hold for the initial data II and III. 
Though all claims in this paper are just numerically verified via three sets of initial data, our methodology can be used to investigate general situations.  All numerical computations in this paper are performed on Mac OS X system  with Intel Core i5 (2.4 GHz, 4 GB) processor.  All initial values and parameters needed in the numerical examples are omitted except in some important examples. 
\end{rem}

\section{The Occurrence of  Two Velocity Impulses} \label{sect_5}

In this section, we look into under what conditions two velocity impulses occur for Problem \ref{problem1} of  free-flight ballistic missiles.
We  say that two velocity impulses occur, which means that two impulse instants are different and meanwhile no one of two velocity impulses is equal to zero.
We first consider the two-impulse space interception problem \ref{problem1} in which the first impulse instant is fixed at $t_0=0$ by using the initial data I.  
\begin{exam}\label{exam0}\rm
Consider Problem \ref{problem1} and the initial data I.
Let the first impulse instant $t_1$ be fixed at zero, i.e., $t_1=t_0=0$. 
Based on the BCs shown in List \ref{table_Nov7_01}, the solver bvp5c gives the solution message
{\small
\begin{verbatim} 
The maximum error is 9.526e-10. Elapsed time is 11689.357989 seconds.
dv1 = [-376.7263; 338.2633; -586.6407], dv2 = 1.0e-10*[-0.3026; 0.2717; -0.4713]
t2 = 3.6313e-10, th = 697.5637
\end{verbatim}
}
\end{exam}

Comparing with the one-impulse solution in Claim \ref{claim2}, we can say that the local optimal two-impulse  solution degenerates to the one-impulse case as $t_1=t_0$  since the order of magnitude of $t_2$ is 1e-10. Hence we raise a question under what conditions a true two-impulse  solution occurs.  We conjecture that if  inequality constraints on impulse instants and the magnitudes of velocity impulses are imposed simultaneously, then two-impulse solutions may occur. It is verified by the following examples. Here all related inequality constraints are converted into equalities by using the static slackness variable method. The corresponding BCs used in this section 
is shown in List \ref{table_Oct10_01} in Appendix C. 

\begin{exam} \label{exam1}\rm 
For the initial data III, we impose inequality constraints component-wise on velocity impulses \eqref{eq_Nov5_01o}
and time constraints on impulse instants 
\begin{align}
\alpha-t_1\leq 0, \quad t_1-\beta\leq 0, \quad \gamma-(t_2-t_1)\leq 0, \quad \alpha,\beta,\gamma\geq 0 \label{eq_Dec23_03}
\end{align}
that is, $t_1\in[\alpha, \beta], t_2-t_1\geq \gamma$.  All parameters and initial values are shown in Table \ref{table06b} in Appendix C.
Then the solver bvp5c gives the solution message 
{\small
\begin{verbatim}
The maximum error is 2.326e-14. Elapsed time is 157.390451 seconds.
dv1 = 1.0e+03*[1.3000; -0.8802; 1.2638]; dv2 = [643.4791; -334.8857; 474.7955] 
abs(dv1)+abs(dv2) = 2882.4177
t_impulse1 = 20, t_impulse2 = 70, t_impact = 3.982766e+02
\end{verbatim}
}
\noindent The resulting position vectors, velocity vectors, their magnitudes, and primer vector and position costate magnitudes are shown in Fig.~\ref{fig_exam3_4_Oct10_04} where the symbols $+, \times, \otimes$ indicate the positions of velocity impulses and  the impact point respectively.
 
\begin{figure}[h!]
 
 \begin{minipage}{8cm}
 \centering
 \includegraphics[width=\hsize]{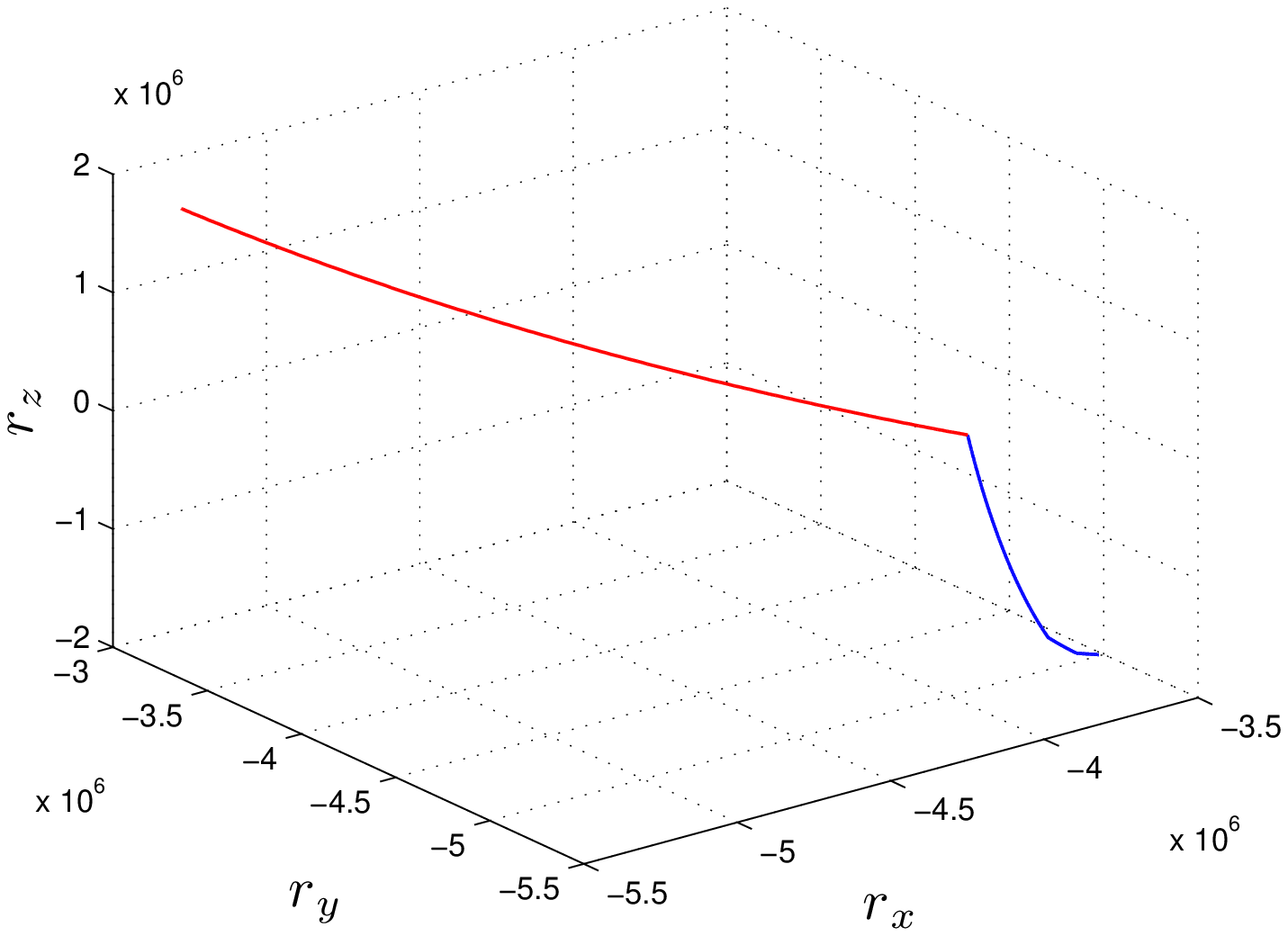} 
\end{minipage}
\qquad 
 \begin{minipage}{8cm}
 \centering
 \includegraphics[width=\hsize]{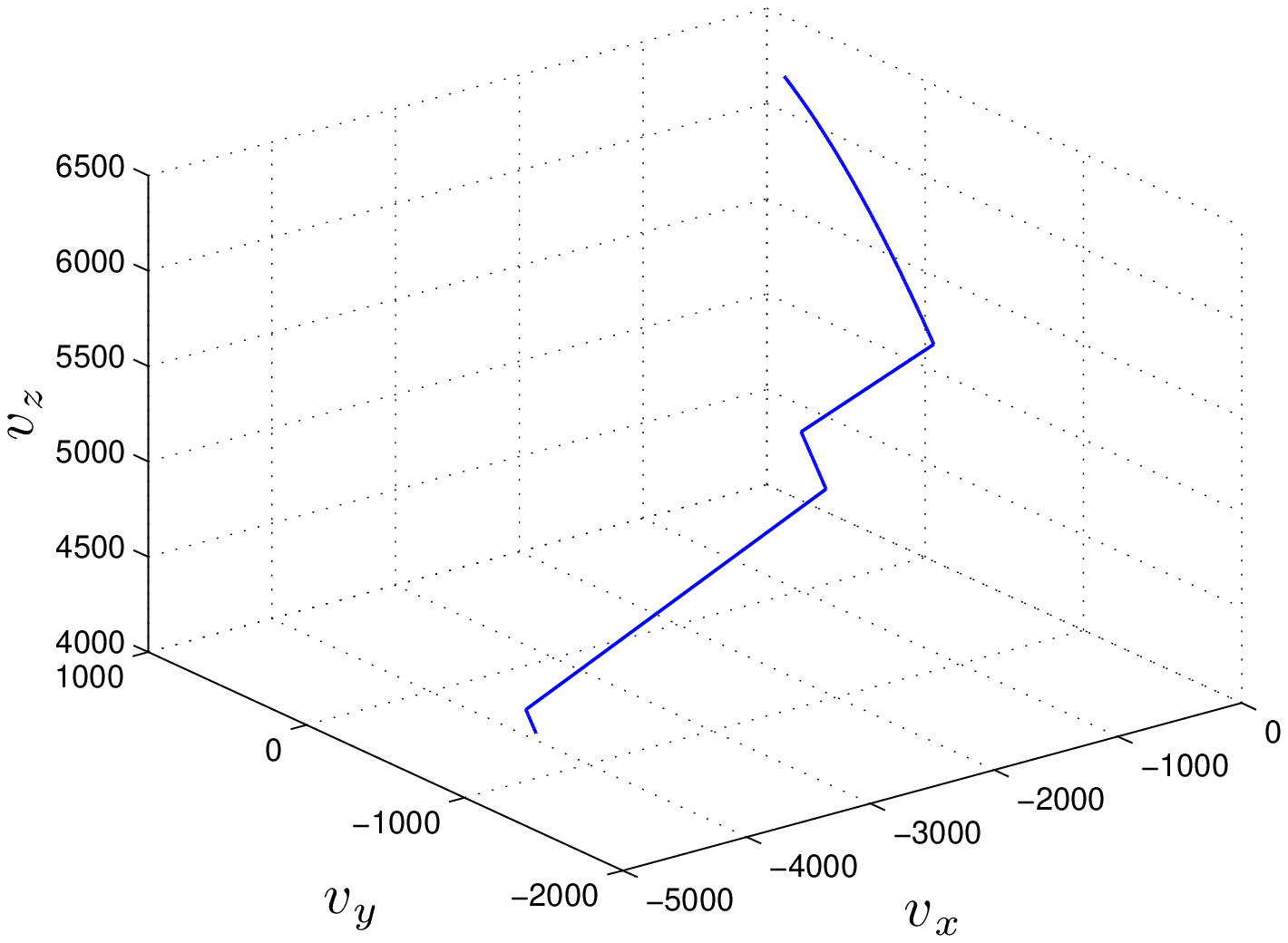}
\end{minipage}

 
 \begin{minipage}{8cm}
 \centering
 \includegraphics[width=\hsize]{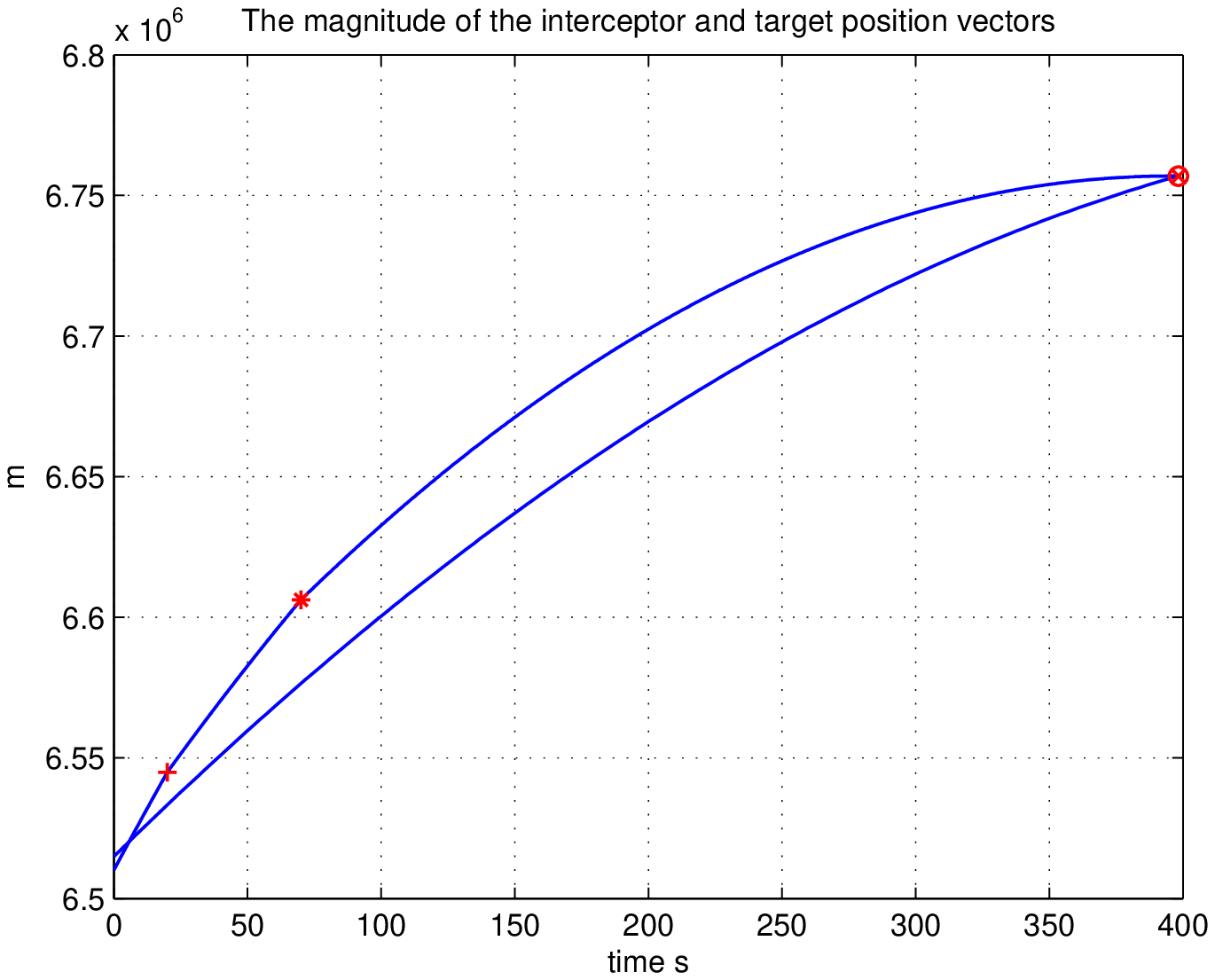} 
\end{minipage}
\qquad 
 \begin{minipage}{8cm}
 \centering
 \includegraphics[width=\hsize]{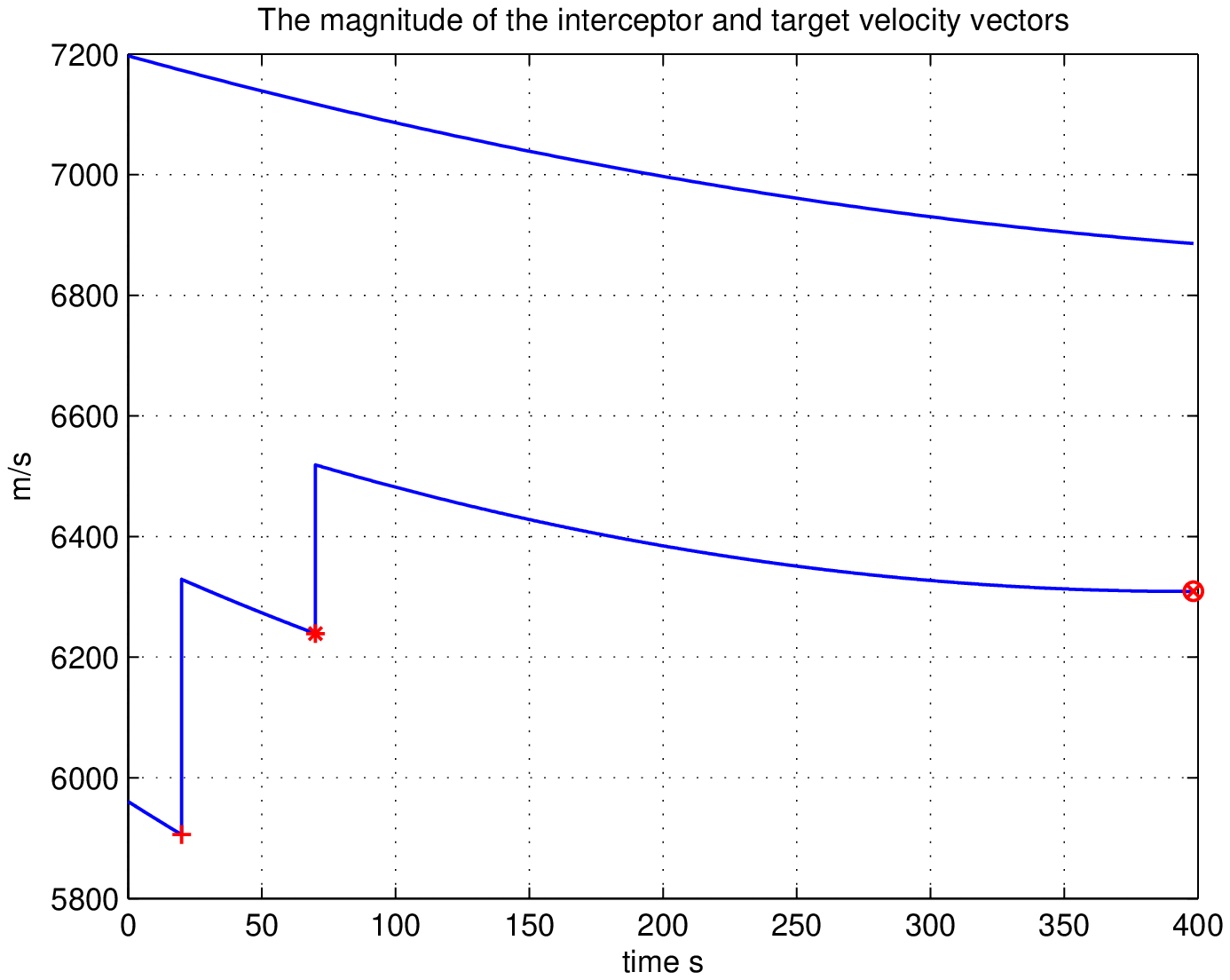}
\end{minipage}

\begin{minipage}{8cm}
 \centering
 \includegraphics[width=\hsize]{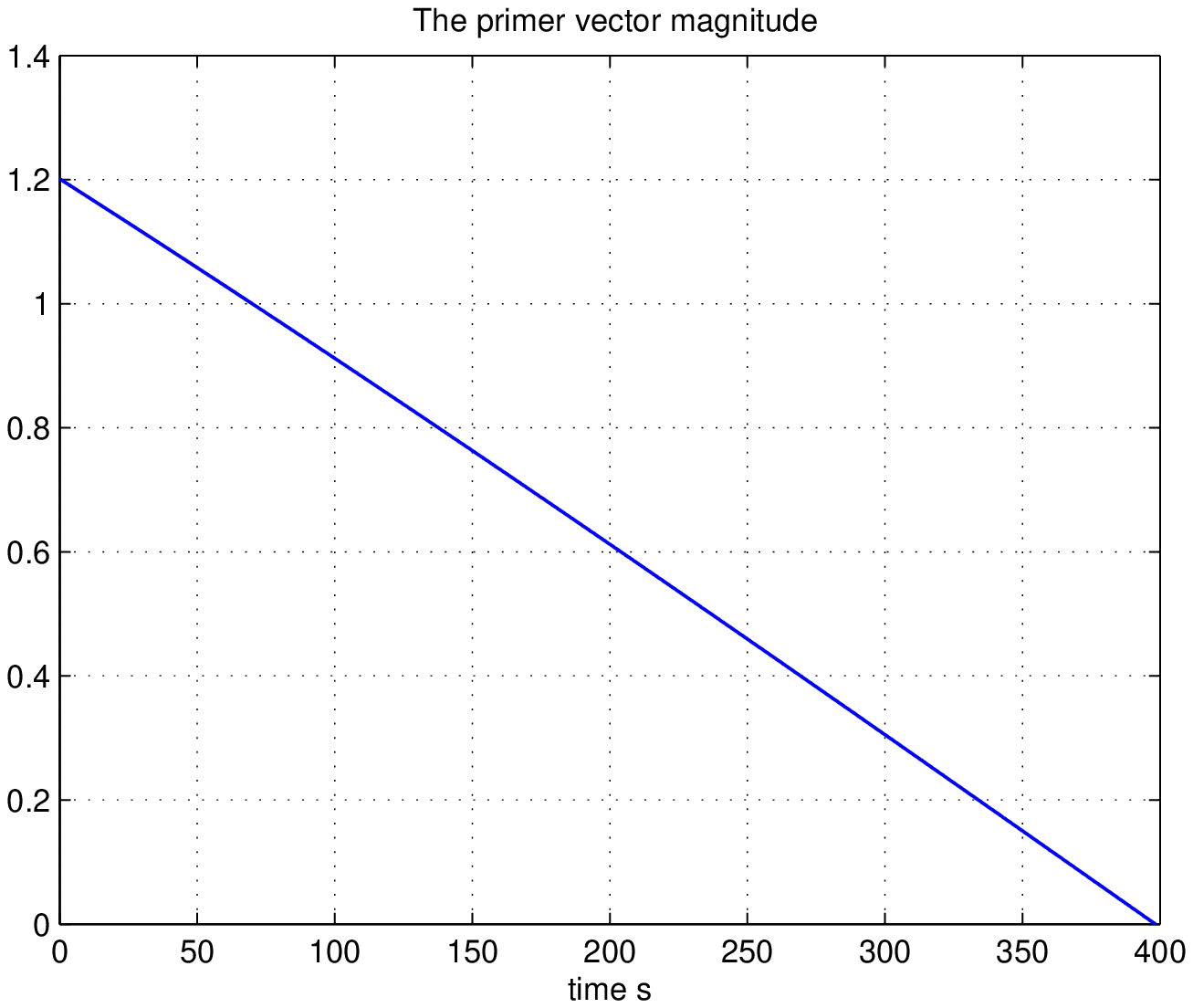} 
\end{minipage}
\qquad 
 \begin{minipage}{8cm}
 \centering
 \includegraphics[width=\hsize]{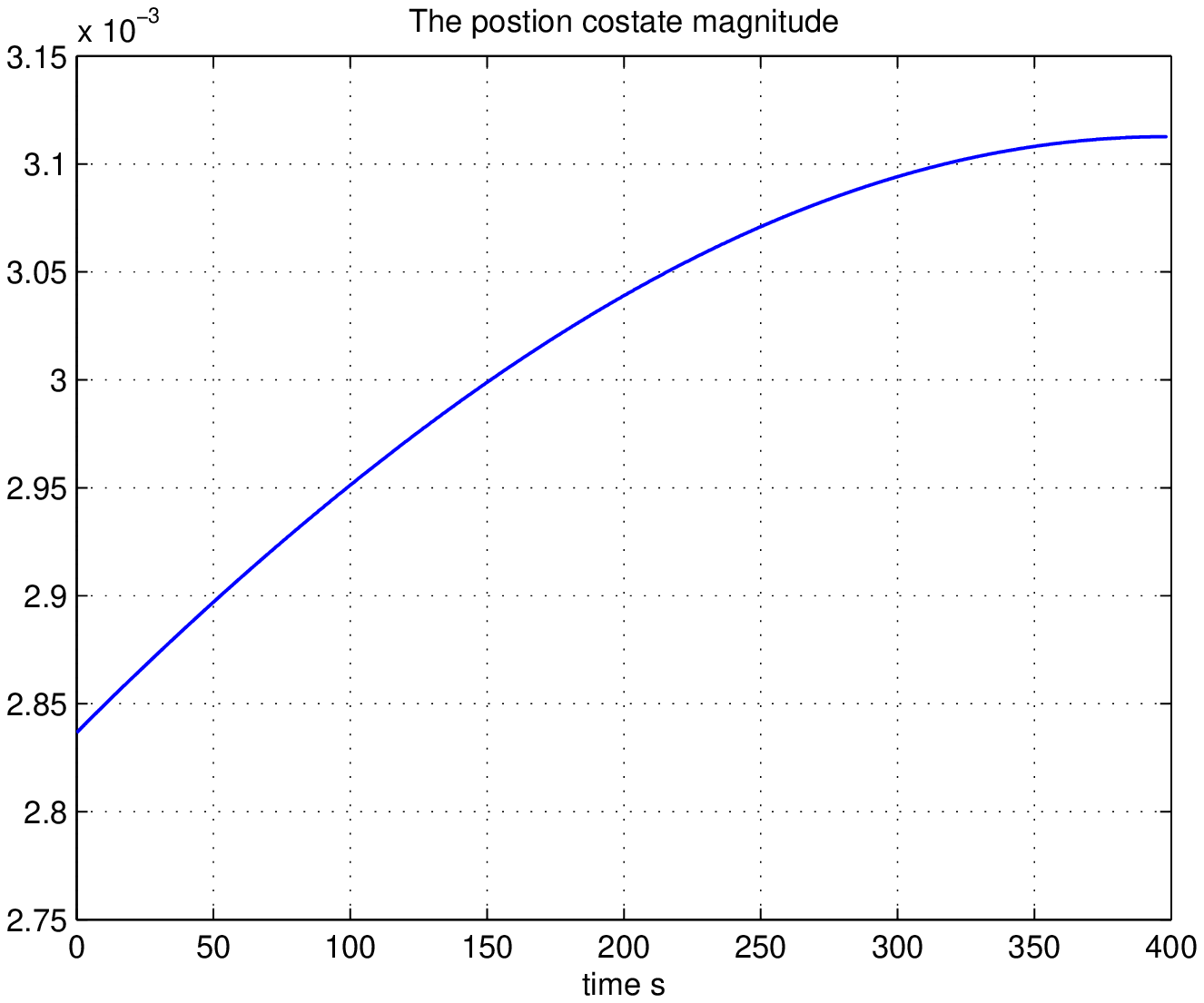}
\end{minipage}
\caption{Position and velocity vectors, their magnitudes, and  primer vector and position costate magnitudes}
\label{fig_exam3_4_Oct10_04}
\end{figure}
\end{exam} 
\begin{exam}\label{exam2}\rm 
For the initial data I, consider the same constraints as in Example \ref{exam1} and set related parameters and initial values by Table \ref{table06a} in Appendix C, then we have 
{\small
\begin{verbatim}
The maximum error is 1.612e-14. Elapsed time is 98.575404 seconds.
dv1 = [-358.9788; 319.8635; -500.0000]; dv2 = [-51.8201; 46.1877; -80.8052]
abs(dv1)+abs(dv2) = 800.1978
t_impulse1 = 20, t_impulse2 = 70, t_impact = 6.957846e+02
\end{verbatim}
}
\end{exam}
Examples \ref{exam1}-\ref{exam2} show that the inequality constraints component-wise on velocity impulses \eqref{eq_Nov5_01o}
and time constraints on impulse instants \eqref{eq_Dec23_03} indeed may result in the occurrence of a local optimal two-impulse solution.  

At the beginning of this section, we have shown that when $t_1=t_0$, a two-impulse case degenerates to the one-impulse one. We use the next example to illustrate that even if there are time and velocity impulse constraints,  this degeneration may still occur. 

\begin{exam} \label{exam3}
\rm Consider the initial data I. We impose the component-wise magnitudes on two velocity impulses as \eqref{eq_Nov5_01m}
and inequality constraints on impulse instants
\begin{align}
0-t_1\leq 0, \quad 0-t_2\leq 0 \label{eq_Dec25_01}
\end{align}
then we have 
{\small
\begin{verbatim}
The maximum error is 1.329e-13. Elapsed time is 107.589378 seconds.
dv1 = [-222.9102; 200.1515; -347.1171], dv2 = [-153.8162; 138.1118; -239.5235]
abs(dv1)+abs(dv2) = 774.9142
t_impulse1=2.411938e-24, t_impulse2=1.030637e-22, t_impact=6.975637e+02
\end{verbatim}
}
\noindent The local optimal impulse instants can be considered as zero since the order of magnitude is 1e-22. Hence 
the local optimal two-impulse solution is actually the solution of the one-impulse case.  

\end{exam}

By Examples \ref{exam1}-\ref{exam3},  it is found that 
\begin{enumerate}
\item the two impulses may occur under time and velocity impulse magnitude constraints; see Examples \ref{exam1}-\ref{exam2}.

\item if the two impulses occur at the instants $t_1$ and $t_2$,  then the magnitude of the two velocity impulse is larger than that of the one-impulse case in which the velocity impulse is applied at $t_1$ instant and less than that of the one-impulse case in which the velocity impulse is applied at $t_2$ instant; see Example \ref{exam2} and Table \ref{table03}. 
 
 \item the upper and lower bounds on velocity impulse magnitudes 
 have effects on the optimal two-impulse solutions. 

 \item with time and velocity impulse magnitude constraints, 
 if the lower bound of two impulse instants is the same, i.e. $t_1,t_2\geq  0$, then the two-impulse case degenerates to the one-impulse one with the constraint $t_1\geq 0$, and also the impulse instant of both cases is equal to $0$; see Example \ref{exam3}. 
 \end{enumerate}
 
The above conclusions are just only drawn from a few numerical examples under consideration, and no theoretical proof is available. Though they seem to be consistent with physical intuitions,  they still need to be verified by more numerical examples. 
To end this section, we make a conjecture about Problem \ref{problem1}.

\begin{conj}\label{con_02} \rm 
For all numerical examples under consideration, the two-impulse space interception problem \ref{problem1} is equivalent to the one-impulse case in the sense of having the same locally optimal solution. 
\end{conj}

We can not obtain a numerical solution to affirm this conjecture within a valid period of time by using the Matlab solver.   
However it is correct if the constraints of velocity impulses are imposed; see Example \ref{exam3}. 

\section{A Terminal Position Constraint} \label{sect_6}

In this section, we consider two-impulse space interception problems with the terminal point constraint, i.e., Problem \ref{problem2} of free-flight ballistic missiles.   
We find that the two-impulse case with free impulse instants is equivalent to the one-impulse case for our examples. Hence we first address the one-impulse case. 
\begin{exam} \label{exam6} \rm
Consider the one-impulse space interception problem with a terminal position constraint for the initial data II. The BCs can be obtained
by removing the BCs related to $t_2$ in List \ref{table_Nov7_02}.
The solution message is given by 
{\small 
\begin{verbatim}
The maximum error is 1.681e-14. Elapsed time is 68.651997 seconds.
dv1 = [-398.4799; 367.8786; -588.7740], abs_dv1 = 800.4847
t_impulse1=5.350987e+01, t_impact=6.824639e+02, tf=9.489139e+02
\end{verbatim}
}

For Problem \ref{problem1} with the initial data II,  the optimal one-impulse solution has the following message
{\small 
\begin{verbatim}
The maximum error is 2.103e-14. Elapsed time is 126.472090 seconds.
dv1 = [-360.3182; 333.9543; -565.8637], abs_dv1 = 749.3707
t_impulse1=4.1561e-27, t_impact=6.835178e+02
\end{verbatim}
}
\begin{figure}[htp!]
 
 \begin{minipage}{8cm}
 \centering
\includegraphics[width=\hsize]{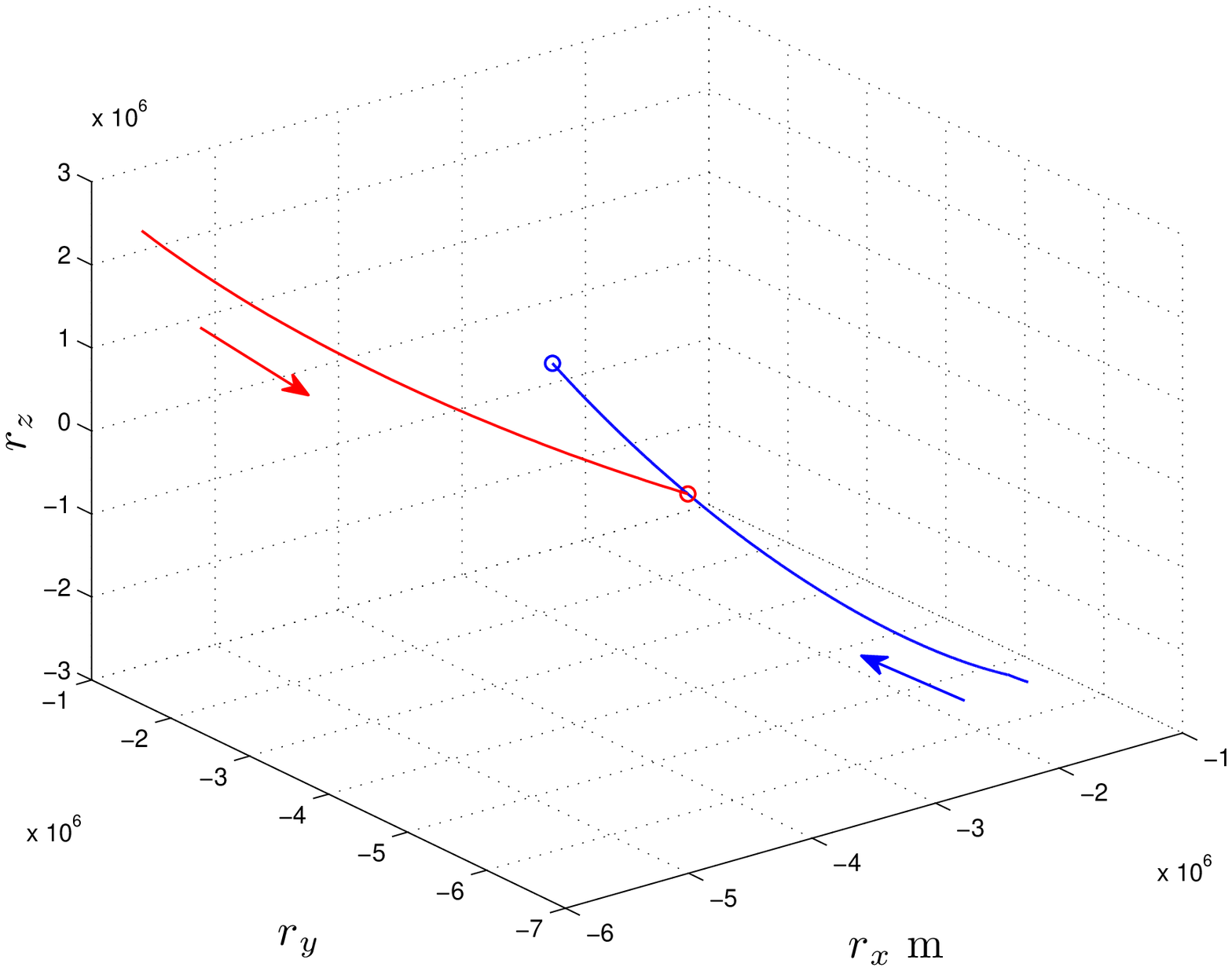}
\end{minipage}
\qquad 
 \begin{minipage}{8cm}
 \centering
\includegraphics[width=\hsize]{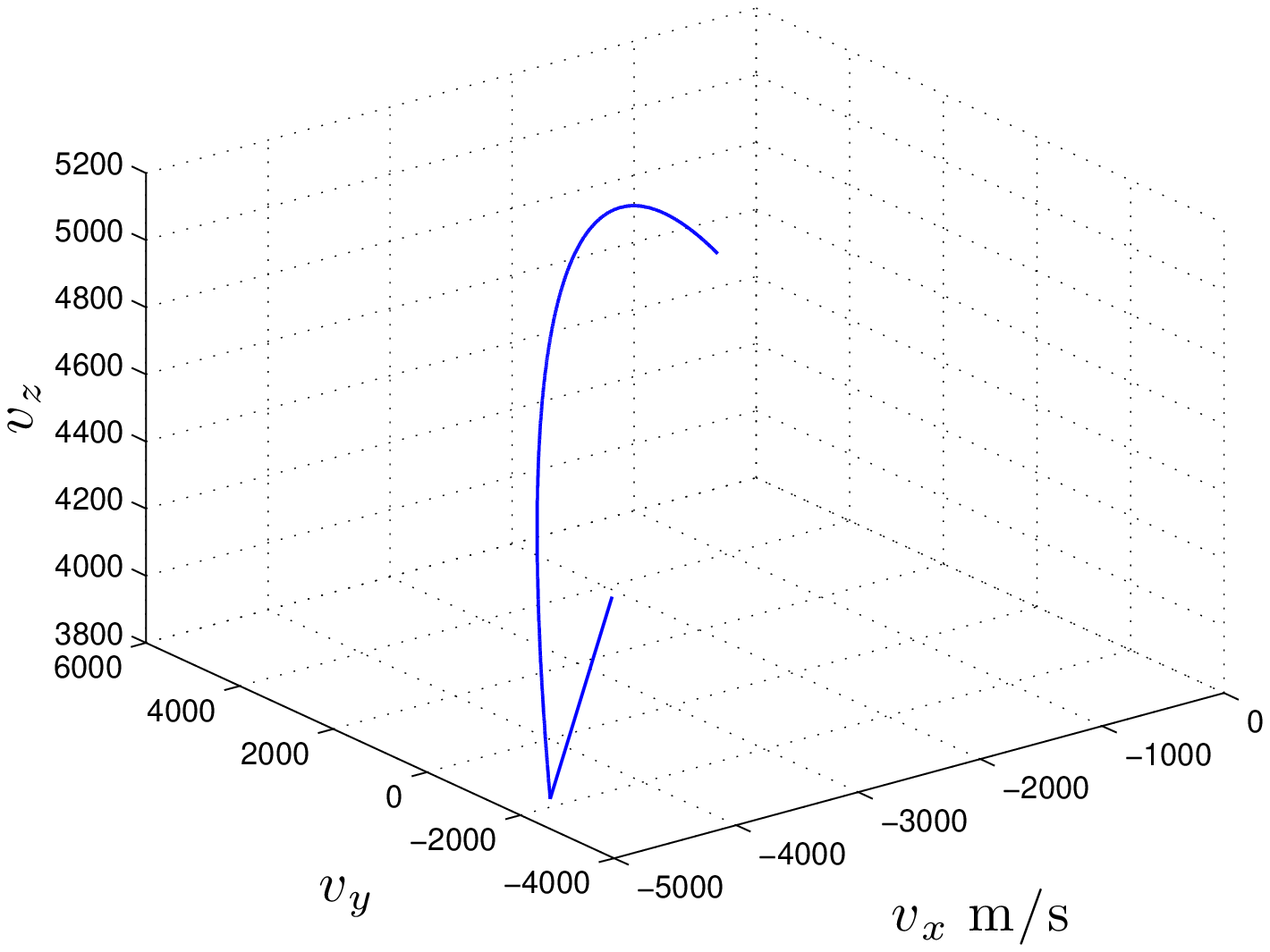}
\end{minipage}
\caption{Interception trajectories and velocity vector}\label{fig_exam3_4_911_l}

\begin{minipage}{8cm}
 \centering
 \includegraphics[width=\hsize]{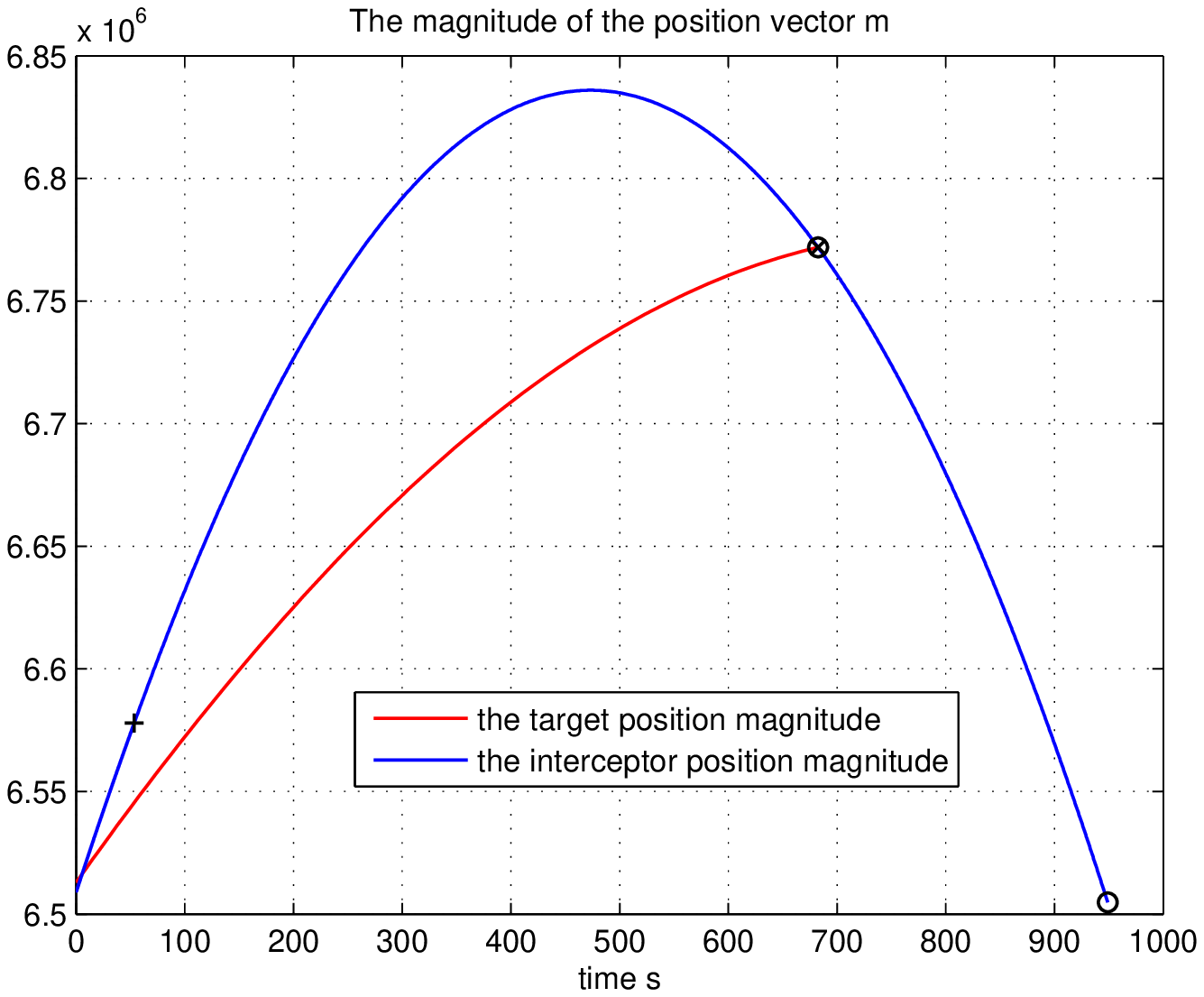}
\end{minipage}
\qquad 
 \begin{minipage}{8cm}
 \centering
  \includegraphics[width=\hsize]{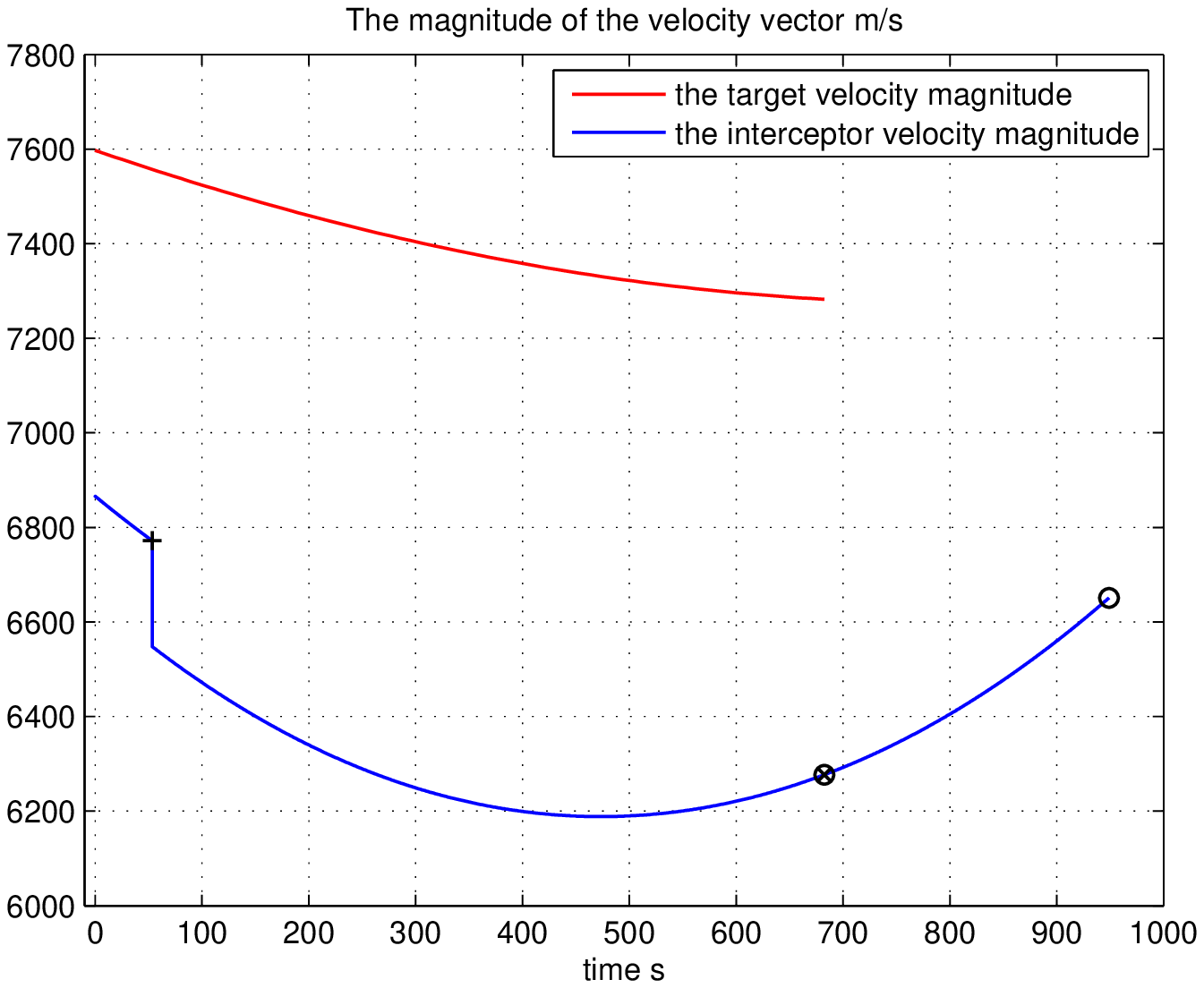}
\end{minipage}
\caption{Position and velocity vector magnitudes}\label{fig_exam3_4_911_m}

 \begin{minipage}{8cm}
 \centering
 \includegraphics[width=\hsize]{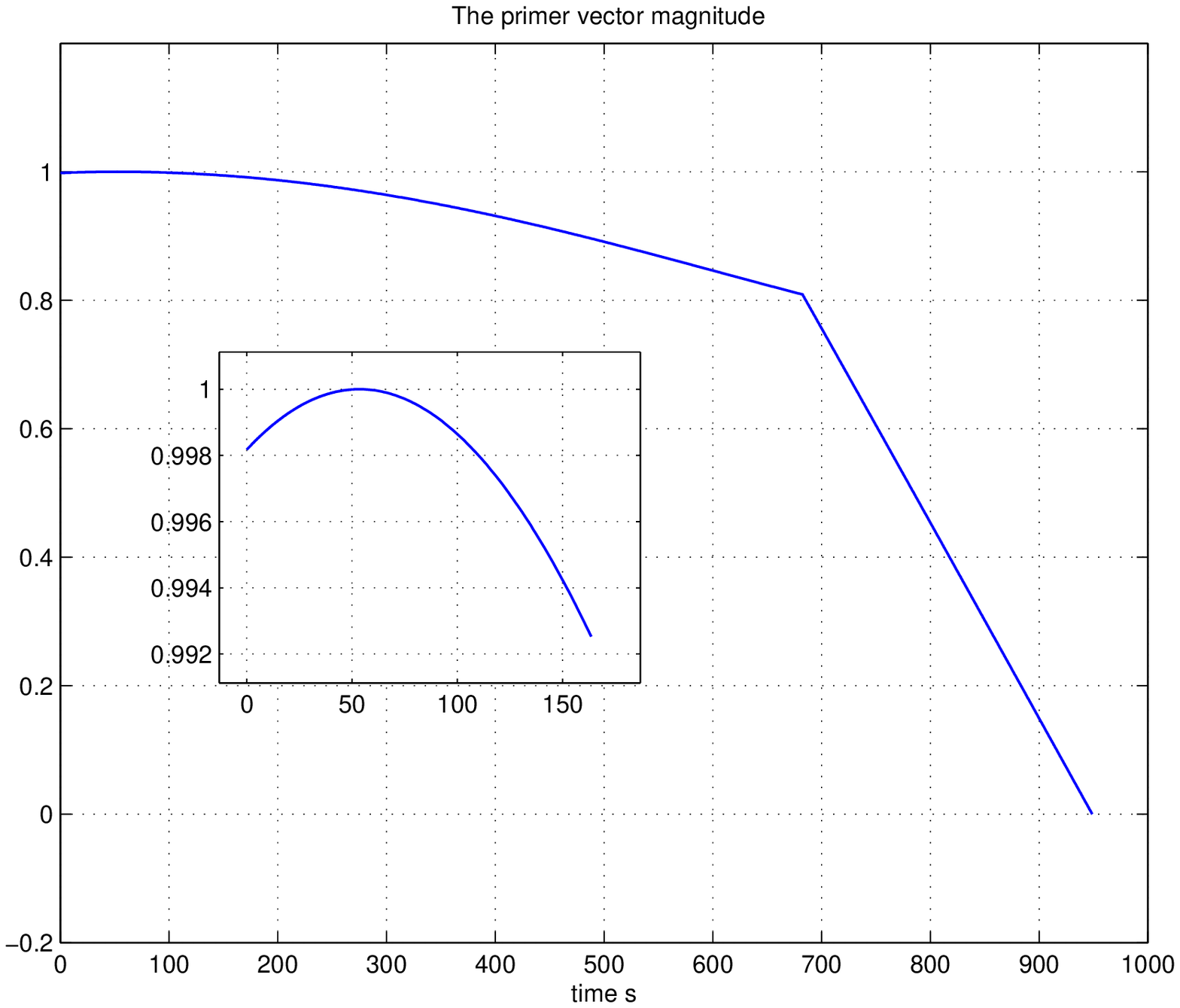}
\end{minipage}
\qquad 
 \begin{minipage}{8cm}
 \centering
  \includegraphics[width=\hsize]{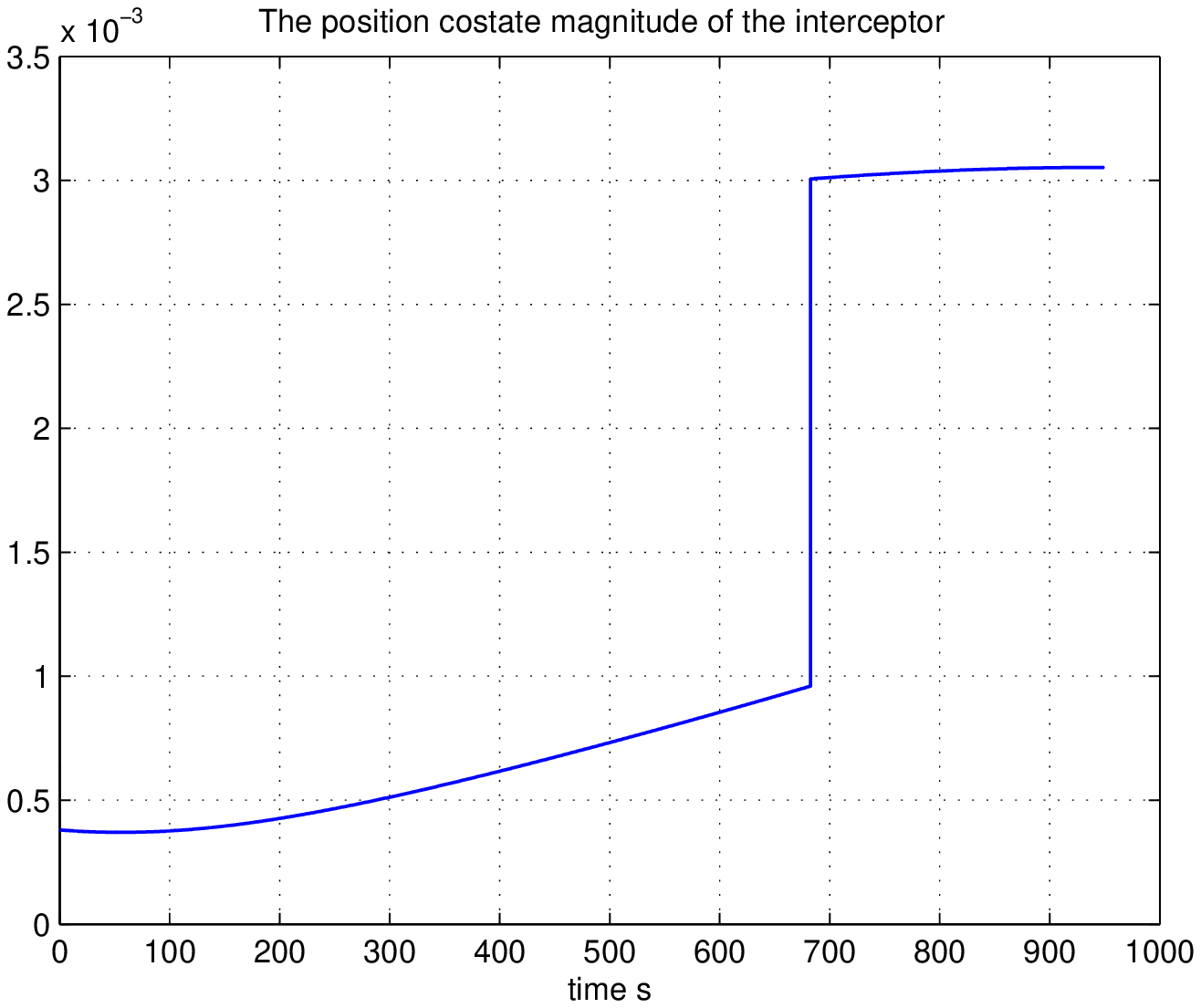}
\end{minipage}
\caption{Primer vector and position costate magnitudes}\label{fig_exam3_4_911_p}
\end{figure}

It is obvious that in this example the cost of Problem \ref{problem1} is less than the cost of Problem \ref{problem2}, and the terminal position constraint indeed changes the impulse instant, the impact instant, the cost, and further the corresponding position and velocity vectors.  Figs.~\ref{fig_exam3_4_911_l}-\ref{fig_exam3_4_911_m} show the trajectories and the velocity vectors of the interceptor and target. 
The primer vector and position costate magnitudes are given in Fig.~\ref{fig_exam3_4_911_p} from which one can see that there is a jump at the impact instant for the position costate and the primer vector achieves its maximal magnitude at the impulse instant $t_1$, where the symbol $\circ$ denotes the terminal position vector (the reference position vector). 
\end{exam}

The next example shows that the two-impulse space interception problem with a terminal position constraint is equivalent to the one-impulse case in Example \ref{exam6}.

\begin{exam} \rm Consider the initial data II for Problem \ref{problem2} of ballistic missiles. The BCs are given by List \ref{table_Nov7_02} and we have the following solution message 
{\small
\begin{verbatim}
The maximum error is 6.284e-15.  Elapsed time is 172.984882 seconds.
dv1 = [-397.8985; 367.3419; -587.9150], dv2 = [-0.5813; 0.5367; -0.8589]
abs(dv1)+abs(dv2): 800.4847
t_impulse1=5.350987e+01, t_impulse2=5.350987e+01, t_impact=6.824639e+02, tf=9.489139e+02
\end{verbatim}
}
It is noted that  $t_2-t_1=7.034373e-13$
which implies the two-impulse case degenerates to the one-impulse case. 
\end{exam}
For Problem \ref{problem2}, it is possible that if $t_1$ is fixed at $t_0$, then an optimal two-impulse solution occurs. 
\begin{exam}\rm
Consider the initial data I. The first impulse instant is fixed at zero. The two-impulse solution is given by
{\small
\begin{verbatim}
The maximum error is 3.308e-10. Elapsed time is 837.853649 seconds.
dv1 = [-358.2908; 321.7015; -558.0683], dv2 = [-18.3931; 16.5315; -28.6650]
t2 = 8.9894e-04, th = 0.7275, tf = 958.9110
abs(dv1)+abs(dv2) = 774.950445
\end{verbatim}
}
The cost is slightly greater than that of the impulse instant free case which equals   774.950428 (Tolerance = 1e-6). 
\end{exam}

\section{Numerical Verification for Multi-constraints} \label{sect_7}


We have known that if time and velocity impulse inequality constraints are imposed, then a local optimal two-impulse solution to Problem \ref{problem1} of  free-flight ballistic missiles may occur. In this section, we consider Problem \ref{problem3} in which there are multi-constraints in terms of inequalities.  The multi-constraints include an additional inequality constraint on the terminal position vector of the interceptor except the time and velocity impulse inequality constraints. 
The static slackness variable method is used to convert the time inequality constraints and the component-wise constraints on velocity impulses into equalities, and the dynamic slackness variable method is used for the terminal  constraints of the interceptor final position vector.
By numerical examples,  we find that the proposed slackness variable methods   successfully solve the two-impulse space interception problem of free-flight ballistic missiles under consideration.  

\begin{exam} \label{exam9}\rm
Consider the two-impulse space interception Problem \ref{problem3} of ballistic missiles with the initial data II. The parameters with respect to all constraints \eqref{eq_Nov5_01m}-\eqref{eq_Nov5_01o} and the initial values of unknown variables are shown in Tables \ref{table10}-\ref{table10c}.
We take the optimal initial value of costates in Problem \ref{problem1} as a guess of costates here
{\small 
\begin{verbatim}
pmr0=1.0e-03*[0.2913 0.2507 -0.0421]; pmv0=[0.4861 -0.4364 0.7571]
\end{verbatim} 
}

\begin{center}\rm
\begin{threeparttable}\small 
\caption{Time constraints \eqref{eq_Nov5_01m} and velocity impulse constraints \eqref{eq_Nov5_01n} }\label{table10}
\begin{tabular}{cccc c ccccc }
\toprule
\multirow{2}{*}{$\alpha$}& \multirow{2}{*}{$\beta$} &\multirow{2}{*}{$\gamma$} &\multirow{2}{*}{$\eta$}  & $ p_{1\min}$ & $ p_{2\min}$  & $ p_{3\min}$  & $ p_{4\min}$  & $ p_{5\min}$  & $ p_{6\min}$  \\
&   &  &  &  $ p_{1\max}$ & $ p_{2\max}$  & $ p_{3\max}$  & $ p_{4\max}$  & $ p_{5\max}$  & $ p_{6\max}$  \\
\hline  
\multirow{2}{*}{$20$} &  \multirow{2}{*}{30} &  \multirow{2}{*}{41} &  \multirow{2}{*}{0}  & -200 &100& -100 & -500 & 100 &-600 \\
&  & & &    100 & 200 & -50 & -100 & 500 &-100
\\
\bottomrule
\end{tabular}
\end{threeparttable}
\end{center}
\begin{center}
\begin{threeparttable}
\caption{Initial $\bm\epsilon, \mathbf{p}_\epsilon$, constraints \eqref{eq_Nov5_01o} and initial $\Delta \mathbf v$} \label{table10b}
\small 
\begin{tabular}{c c ccccc }
\toprule
\multirow{2}{*}{$\epsilon_x, \epsilon_y , \epsilon_z, p_{x}, p_y, p_z$} & $ r_{x\min}$ & $ r_{y\min}$  & $ r_{z\min}$  & $\text{dv1}_{x}$  & $ \text{dv1}_{y}$  & $ \text{dv1}_{z}$  \\
&  $ r_{x\max}$ & $ r_{y\max}$  & $ r_{z\max}$  & $ \text{dv2}_{x}$  & $ \text{dv2}_{y}$  & $\text{dv2}_{z}$  \\
\hline
\multirow{2}{*}{1}  & -500 &-500& -500 & -100 & 100 &-100 \\
&      500 & 500 & 500 & -300 & 300 &-400
\\
\bottomrule
\end{tabular}
\end{threeparttable}
\end{center}
\begin{center}\rm
\begin{threeparttable}
\caption{Initial scaled instants, weighting coefficients in \eqref{eq_Dec15_01} and Lagrange multipliers}\label{table10c}
\small 
\begin{tabular}{cccc c c c c}
\toprule
$t_1$  &  $t_2$ &  {$t_h$} &  {$t_f$} &  {$\lambda_1; \lambda_2, \lambda_3, \lambda_4$} &  $\mu_1, \mu_2, \dots,\mu_{12}$ &  {$k_{3x}, k_{3y}, k_{3z},  k_{4x}, k_{4y}, k_{4z}$}& $\eta_1, \eta_2, \dots,\eta_{6}$ \\
\hline
0.02 & 0.04   & 0.14 & 950 & 1.006; 1 &  {1} & {1} & -2
\\ 
\bottomrule
\end{tabular}
\end{threeparttable}
\end{center}
There are 116 BCs, the part of which related to the costates is given by List \ref{table_Oct21_01}. 
Then by trial and error, the Matlab solver {\rm \texttt{bvp5c}} provides successfully with the message of the solution to Problem \ref{problem3} in which we set Tolerance = 1e-9
 {\small
\begin{verbatim}
The maximum error is 9.887e-10. Elapsed time is 3841.473078 seconds. 
dv1 = [-100; 100; -100], dv2 = [-297.2786; 265.4892; -490.6278]
abs(dv1)+abs(dv2) = 805.3242
t_impulse1=2.075471e+01, t_impulse2=6.175471e+01, t_impact=6.824994e+02, tf=9.491821e+02
The difference between the optimal position at the final time and the reference position:
500.0000, 500.0000, 500.0000
\end{verbatim}
}

\begin{figure}[h]
 \begin{minipage}{8cm}
 \centering
 \includegraphics[width=\hsize]{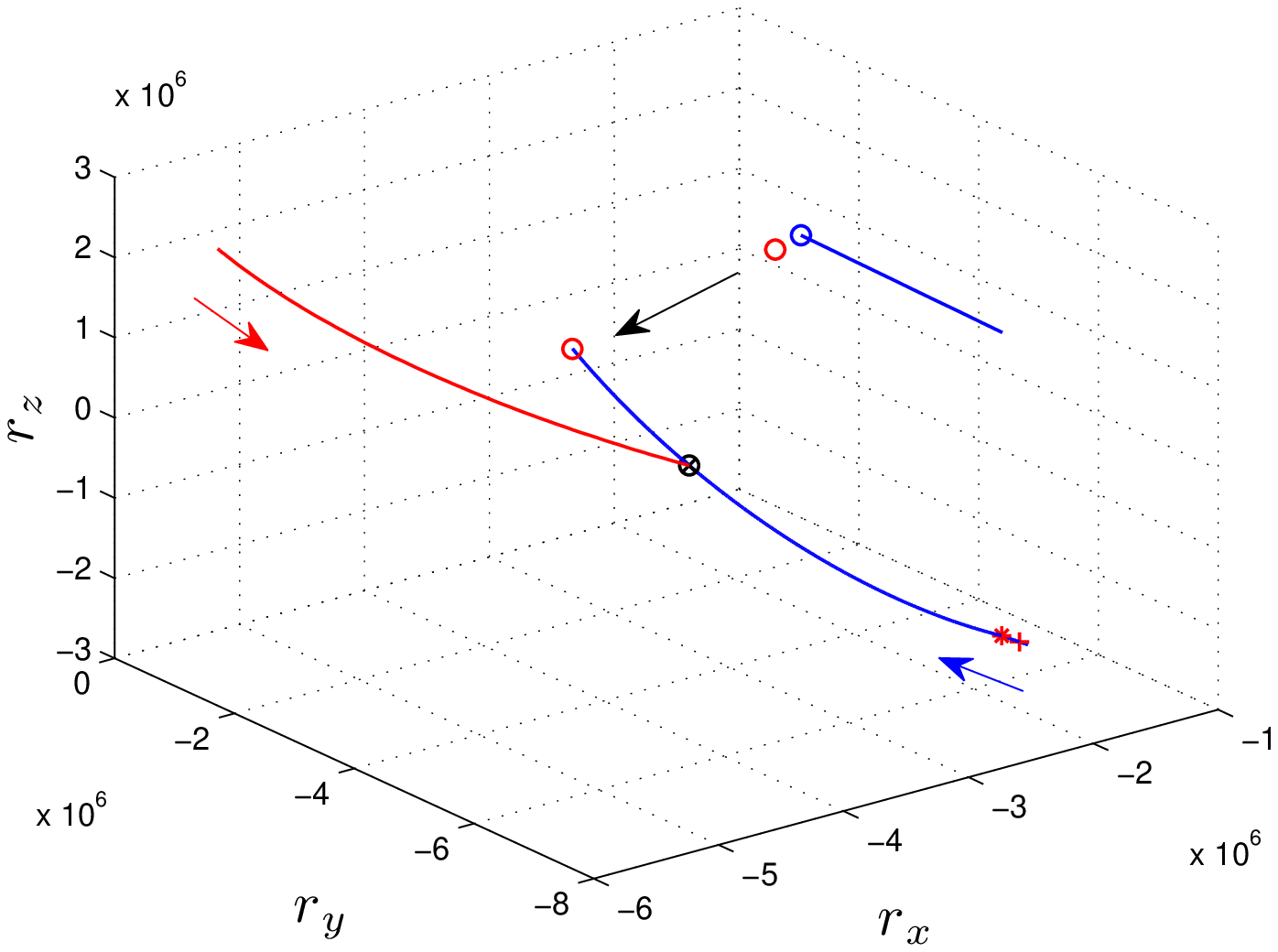} 
\end{minipage}
\quad 
 \begin{minipage}{8cm}
 \centering
 \includegraphics[width=\hsize]{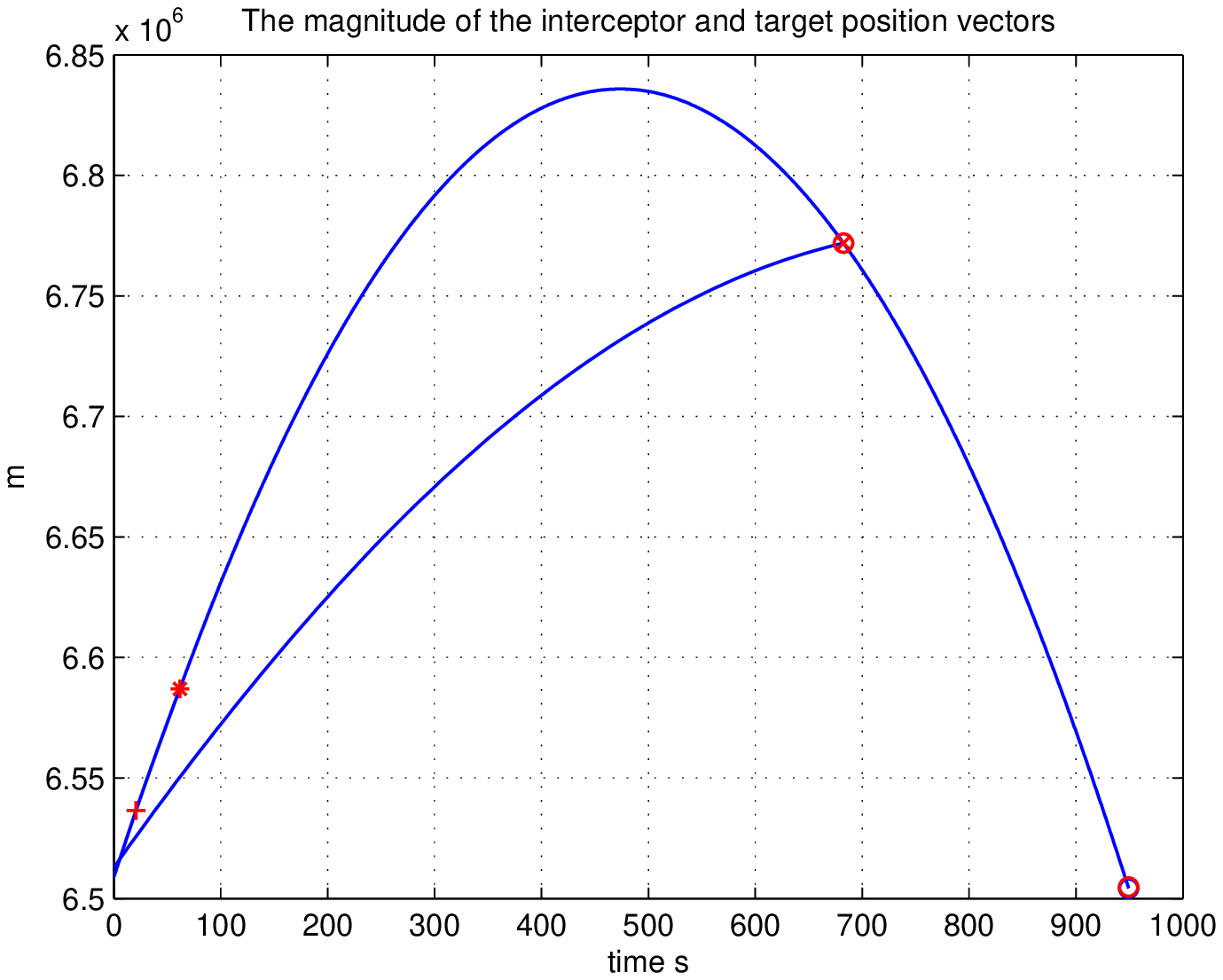}
\end{minipage}
 \begin{minipage}{8cm}
 \centering
 \includegraphics[width=\hsize]{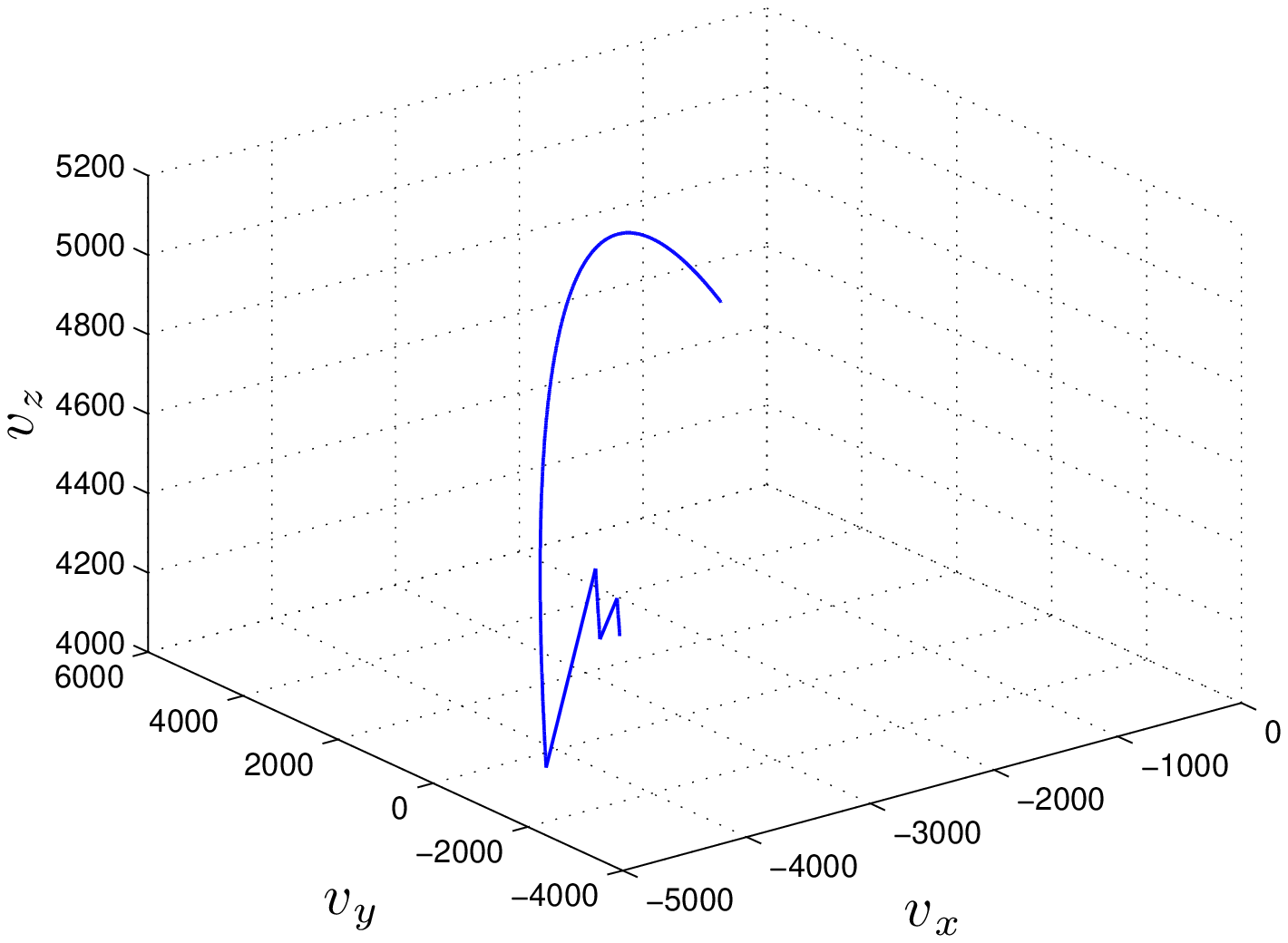} 
\end{minipage}
\quad 
 \begin{minipage}{8cm}
 \centering
 \includegraphics[width=\hsize]{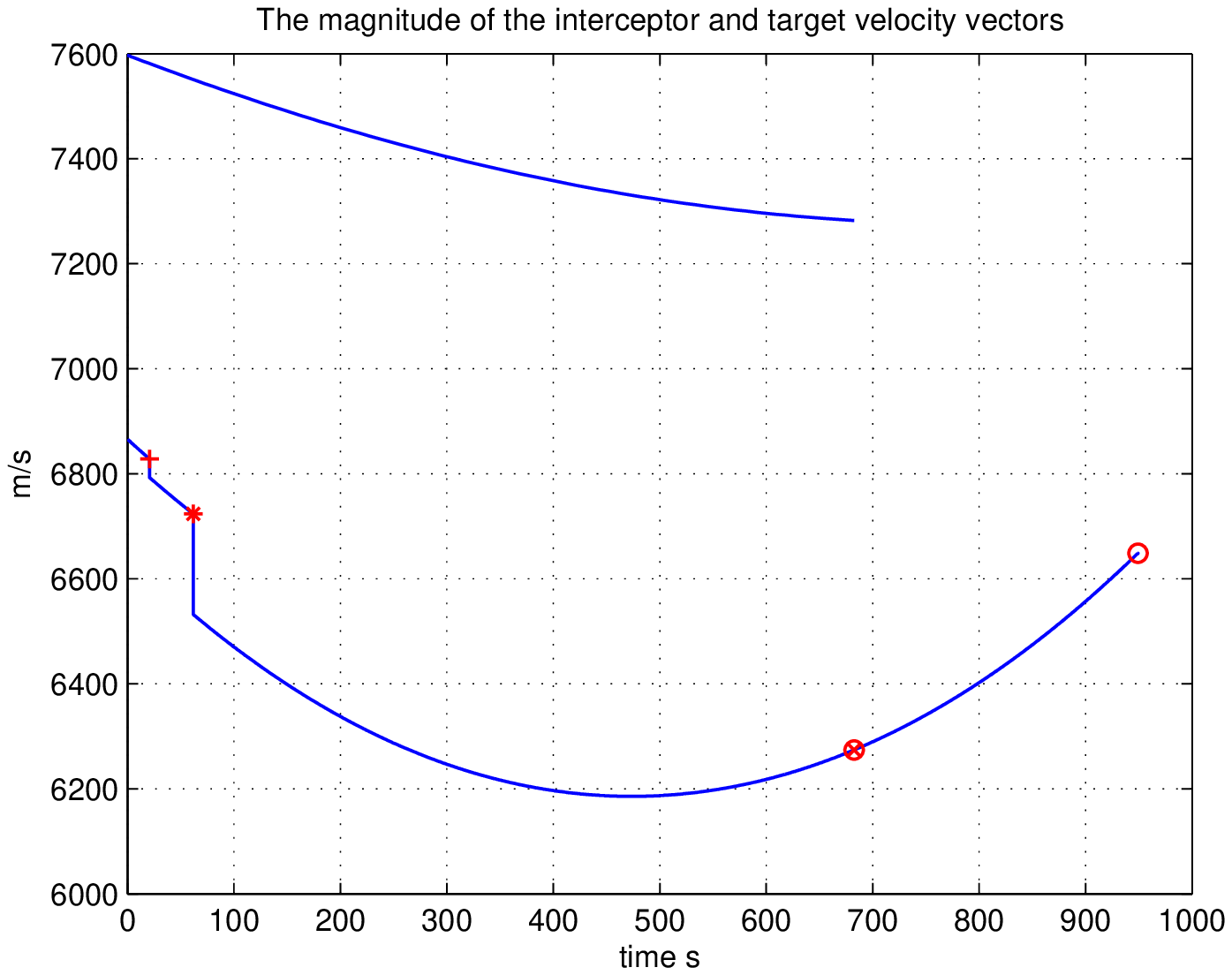}
\end{minipage}
\caption{Position and velocity vectors and their magnitudes}
\label{fig_exam3_4_Oct23_04f}
\end{figure}

\begin{figure}[h]
 \begin{minipage}{8cm}
 \centering
\includegraphics[width=\hsize]{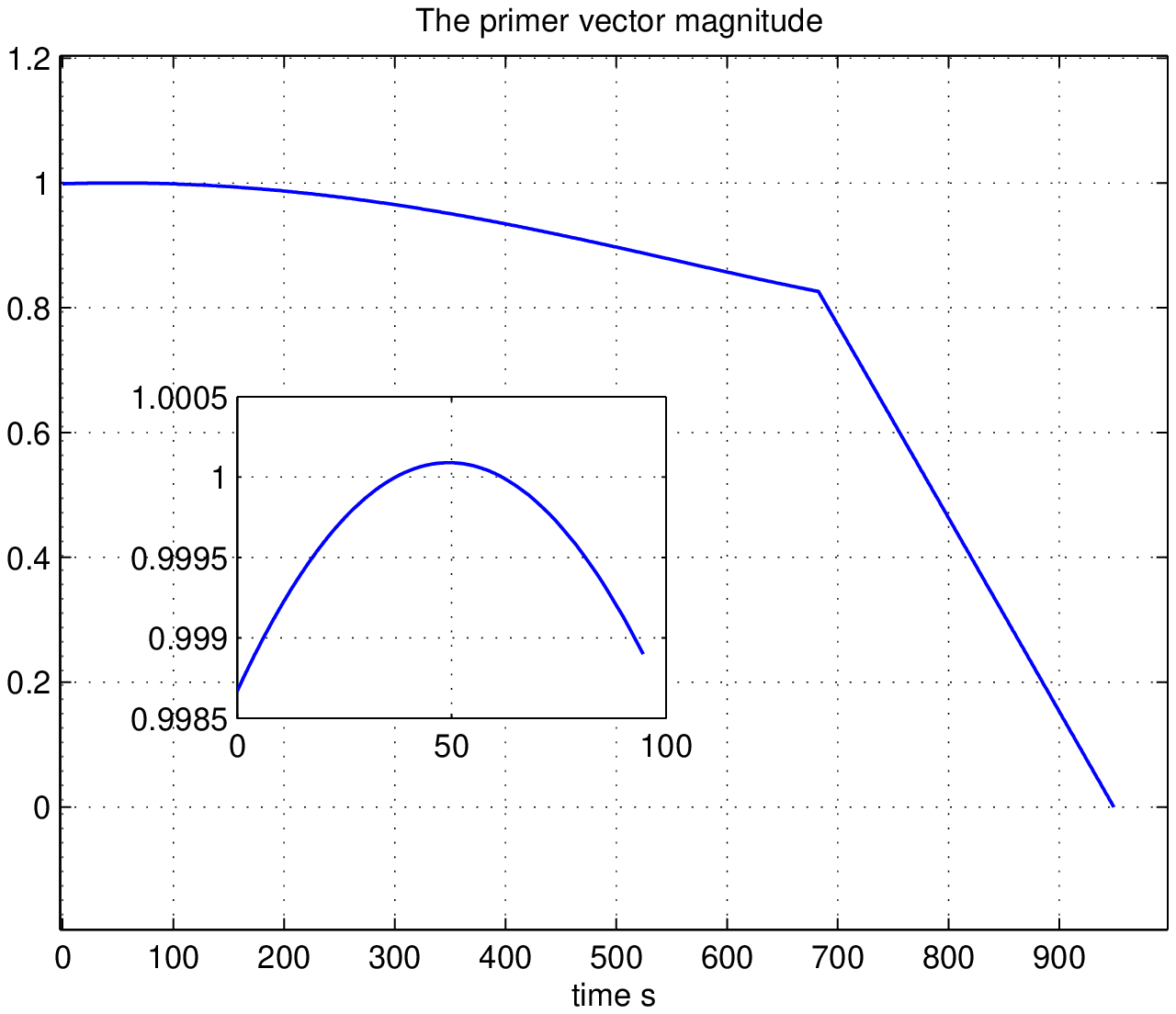} 
\end{minipage}
\quad 
 \begin{minipage}{8cm}
 \centering
 \includegraphics[width=\hsize]{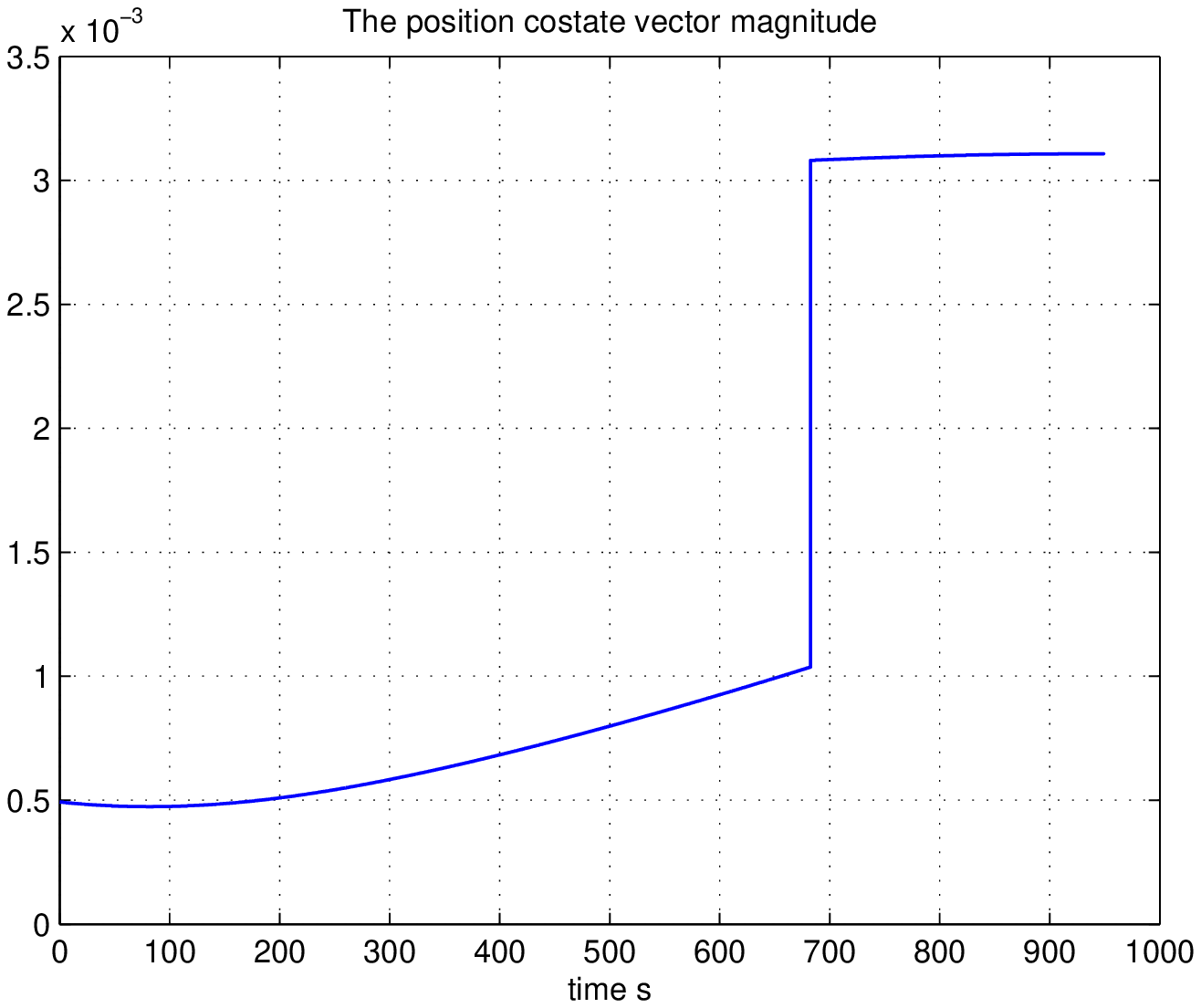}
\end{minipage}

\begin{minipage}{8cm}
 \centering
 \includegraphics[width=\hsize]{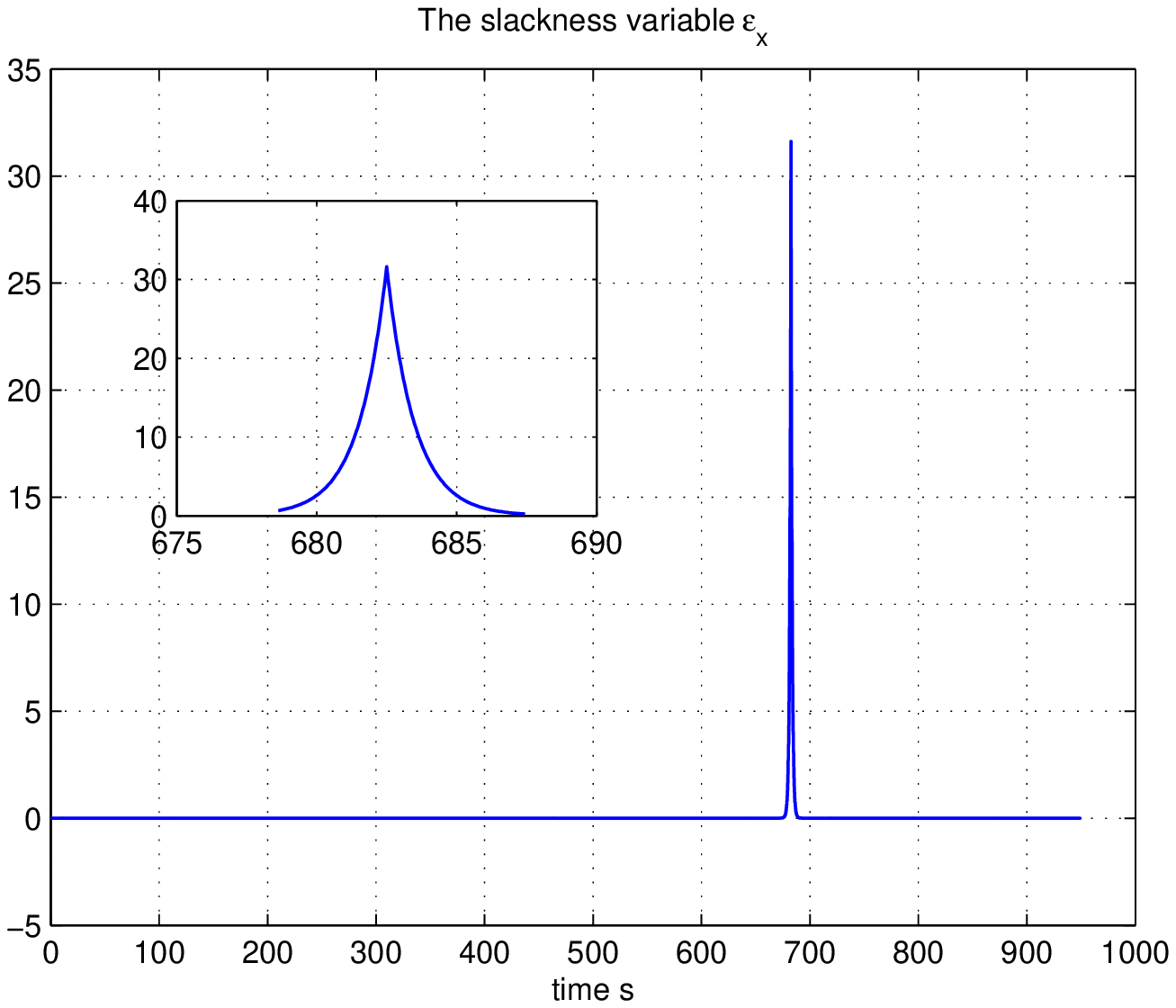} 
\end{minipage}
\qquad 
 \begin{minipage}{8cm}
 \centering
 \includegraphics[width=\hsize]{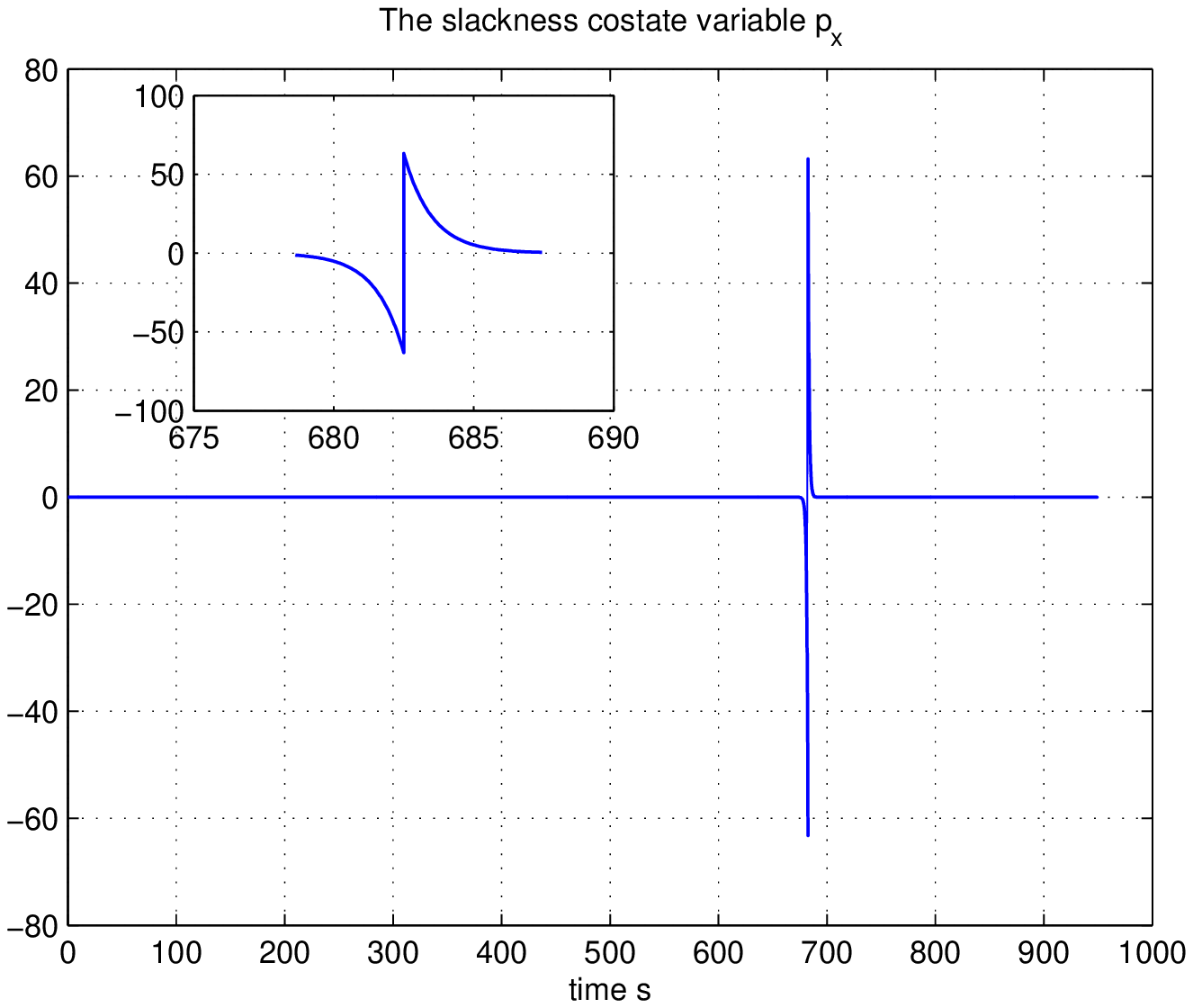}
\end{minipage}

\caption{Magnitudes of primer vector and position costate, and slackness variable $\epsilon_x$ and its costate $p_{x\epsilon}$ }
\label{fig_exam3_4_Oct23_04}
\end{figure}

 Figs.~\ref{fig_exam3_4_Oct23_04f}-\ref{fig_exam3_4_Oct23_04} show the position vector, the velocity vector and their magnitudes for both interceptor and target, and the primer vector and position costate magnitudes of the interceptor, respectively.  The symbols $*, +, \circ, \otimes $  correspond to two impulse instants, the terminal point of the interceptor, and the impact instant, respectively.  If we use the Matlab solver bvp4c instead, assume its tolerance to be 1e-6, set $\lambda_1=1.0056$,  and modify the code by vectorization methods, then a solution approximately equal to the one shown above is found satisfying all constraints and the computation time is about
70 seconds.

\end{exam}

It is found that the two-impulse space interception problem with multi-constraints cannot be solved within a valid period of time if the static slackness variable method  is used for the terminal position inequality constraints.  Hence instead, the dynamic slackness variable method is introduced in this paper. One can see from Example \ref{exam9} that the solution given by the dynamic slackness variable method satisfies all constraints. 

\section{Conclusion}\label{Sect_8}

Using the calculus of variations, we have solved two-impulse space interception problems with multi-constraints proposed in this paper. 
The multi-constraints include impulse and impact instant constraints, terminal  constraints on the final position of the interceptor, and component-wise constraints on the magnitudes of velocity impulses.
A number of conclusions concerning two-impulse space interception problems have been established based on highly accurate numerical solutions.

One question confuses us when we deal with two-impulse interception problems of ballistic missiles by using Matlab boundary value problem solvers: Under what circumstance a true two-impulse optimal solution occurs for Problem \ref{problem1}?
By numerical solutions, we first find that if the first impulse instant is fixed at the initial time and the second one is free, then its optimal solution degenerates to the one-impulse case in which the optimal pulse instant is equal to the initial time, which means that the one-impulse solution is optimal. 
We then conjecture that to ensure that a two-impulse optimal solution occurs, we must impose constraints on impulse instants and the magnitudes of velocity impulses.  Numerical solutions affirm this conjecture.  
Therefore we propose Problem \ref{problem3}. Then the static and dynamic slackness variable methods are introduced in order to solve Problem \ref{problem3} by using Matlab boundary value problem solvers.  Numerical results show that our method makes senses for the two-impulse space interception problems of ballistic missiles with multi-constraints.

From the perspective of optimal control theory, numerical methods are categorized into two approaches to solving optimization problems: indirect and direct methods. 
In an indirect method, an optimal solution is obtained by satisfying a set of necessary conditions. 
Conversely, a direct method just only uses its cost to obtain an optimal solution \cite{prussing_2010}.  
Here we use the indirect method to numerically solve optimal space interception problems. The advantage of an indirect method is its high solution accuracy. The shortage is also obvious, that is, 
almost all examples are solved by trial and error.  
In addition, Conjecture \ref{con_02} about the equivalence of Problem 1 and the one-impulse case still remains open. It is correct if the constraints of velocity impulses are imposed. 

\section*{Appendix A: The derivation of the BCs for Problem \ref{problem1}} \label{App_1}

The first variation of the augmented cost functional
\begin{align} \label{eq_Dec5_01}
\delta \tilde J&=\dfrac{\Delta\mathbf{v}_0^{\text T} }{\abs{\Delta \mathbf{v}_0}}\delta \mathbf{v}_0+
\dfrac{\Delta\mathbf{v}_1^{\text T} }{\abs{\Delta \mathbf{v}_1}}\delta \mathbf{v}_1\notag\\
&\quad +d\mathbf{r}_T^{\text T}(t_h)\dfrac{\partial g_{h} (\mathbf{r}_T(t_h))}{\partial \mathbf{r}_T(t_h)}\Big |_*\gamma_{h}
+d \mathbf{r}_M^{\text T}(t_h)\dfrac{\partial g_{h} (\mathbf{r}_M(t_h), \mathbf{r}_T(t_h))}{\partial \mathbf{r}_M(t_h)}\Big |_*\gamma_{h}
\notag\\
&\quad+\mathbf{q}_{v1}^{\text T}\left [d \mathbf v(t_0^{+})-d \mathbf v(t_0^{-})-\delta \mathbf{v}_0\right]
+\mathbf{q}_{v2}^{\text T}\left[d \mathbf v(t_1^{+})-d \mathbf v(t_1^{-})-\delta \mathbf{v}_1\right ]
\notag\\
&(*) 
\begin{cases}
+\left (H_{M1*}-\mathbf p_{Mr1}^{\text T} \dot{\mathbf r}_M-\mathbf p_{Mv1}^{\text T} \dot{\mathbf v}_M\right)\Big |_{t_1^{-*}}\delta t_1
-\left (H_{M2*}-\mathbf p_{Mr2}^{\text T} \dot{\mathbf r}_M-\mathbf p_{Mv2}^{\text T} \dot{\mathbf v}_M\right)\Big |_{t_1^{+*}}\delta t_1\notag\\
+
\left (H_{M2*}-\mathbf p_{Mr2}^{\text T} \dot{\mathbf r}_M-\mathbf p_{Mv2}^{\text T} \dot{\mathbf v}_M\right)\Big |_{t_h^{-*}}\delta t_h\notag\\
\end{cases}\\
&(\dagger)
\begin{cases}
-\mathbf p_{Mr1}^{\text T}(t_1^{-*})\delta \mathbf r_M(t_1^{-*})+\mathbf p_{Mr2}^{\text T}(t_1^+)\delta \mathbf r_M(t_1^{+*})
-\mathbf p_{Mr2}^{\text T}(t_h^{-*})\delta \mathbf r_M(t_h^{-*})\notag\\
+\mathbf p_{Mv1}^{\text T}(t_0^{+})\delta \mathbf v_M(t_0^{+})
 -\mathbf p_{Mv1}^{\text T}(t_1^{-*})\delta \mathbf v_M(t_1^{-*})+\mathbf p_{Mv2}^{\text T}(t_1^{+*})\delta \mathbf v_M(t_1^{+*})-\mathbf p_{Mv2}^{\text T}(t_h^{-*})\delta \mathbf v_M(t_h^{-*})
\notag\\
\end{cases}\notag
\\
&(**) 
\begin{cases}
+\left (H_{T1*}-\mathbf p_{Tr1}^{\text T} \dot{\mathbf r}_T-\mathbf p_{Tv1}^{\text T} \dot{\mathbf v}_T\right)\Big |_{t_1^{-*}}\delta t_1\notag\\
-\left (H_{T2*}-\mathbf p_{Tr2}^{\text T} \dot{\mathbf r}_T-\mathbf p_{Tv2}^{\text T} \dot{\mathbf v}_T\right)\Big |_{t_1^{+*}}\delta t_1
+
\left (H_{T2*}-\mathbf p_{Tr2}^{\text T} \dot{\mathbf r}_T-\mathbf p_{Tv2}^{\text T} \dot{\mathbf v}_T\right)\Big |_{t_h^{-*}}\delta t_h
\notag\\
\end{cases}
\notag\\
&(\dagger\dagger)
\begin{cases}
-\mathbf p_{Tr1}^{\text T}(t_1^{-*})\delta \mathbf r_T(t_1^{-*})+\mathbf p_{Tr2}^{\text T}(t_1^+)\delta \mathbf r_T(t_1^{+*})-\mathbf p_{Tr2}^{\text T}(t_h^{-*})\delta \mathbf r_T(t_h^{-*})\\
 -\mathbf p_{Tv1}^{\text T}(t_1^{-*})\delta \mathbf v_T(t_1^{-*})+\mathbf p_{Tv2}^{\text T}(t_1^{+*})\delta \mathbf v_T(t_1^{+*})-\mathbf p_{Tv2}^{\text T}(t_h^{-*})\delta \mathbf v_T(t_h^{-*})
\notag
\end{cases}
\\
&\quad +\int_{t_0^+}^{t_1^{-*}}{\Bigg [}\left (\dfrac{\partial H_{M1}(\mathbf r_M,\mathbf v_M, \mathbf p_{M1})}{\partial \mathbf r_M}+\dot{\mathbf p}_{Mr1}(t)\right )^{\text T}_*\delta \mathbf r_M
+\left(\dfrac{\partial H_{M1}(\mathbf r_M,\mathbf v_M, \mathbf p_{M1})}{\partial \mathbf v_M}+\dot{\mathbf p}_{Mv1}(t)\right)^{\text T}_*\delta \mathbf v_M{\Bigg ]}\text{d} t\notag\\
&\quad +\int_{t_1^+}^{t_h^{-*}}\left [\left (\dfrac{\partial H_{M2}(\mathbf r_M,\mathbf v_M, \mathbf p_{M2})}{\partial \mathbf r_M}+\dot{\mathbf p}_{Mr2}(t)\right )^{\text T}_*\delta \mathbf r_M+\left(\dfrac{\partial H_{M2}(\mathbf r_M,\mathbf v_M, \mathbf p_{M2})}{\partial \mathbf v_M}+\dot{\mathbf p}_{Mv2}(t)\right)^{\text T}_*\delta \mathbf v_M\right ]\text{d} t\notag\\
&\quad +\int_{t_0^+}^{t_1^{-*}}\left [\left (\dfrac{\partial H_{T1}(\mathbf r_T,\mathbf v_T, \mathbf p_{T1})}{\partial \mathbf r_T}+\dot{\mathbf p}_{Tr1}(t)\right )^{\text T}_*\delta \mathbf r_T+\left(\dfrac{\partial H_{T1}(\mathbf r_T,\mathbf v_T, \mathbf p_{T1})}{\partial \mathbf v_T}+\dot{\mathbf p}_{Tv1}(t)\right)^{\text T}_*\delta \mathbf v_T\right ]\text{d} t\notag\\
&\quad +\int_{t_1^+}^{t_h^{-*}}\left [\left (\dfrac{\partial H_{T2}(\mathbf r_T,\mathbf v_T, \mathbf p_{T2})}{\partial \mathbf r}+\dot{\mathbf p}_{Tr2}(t)\right )^{\text T}_*\delta \mathbf r_T+\left(\dfrac{\partial H_{T2}(\mathbf r_T,\mathbf v_T, \mathbf p_{T2})}{\partial \mathbf v_T}+\dot{\mathbf p}_{Tv2}(t)\right)^{\text T}_*\delta \mathbf v_T\right ]\text{d} t 
\end{align}
where we use $d(\cdot)$ to denote the difference between  
the varied path and the optimal path taking into account the differential change in a time instant, 
e.g.,
\[
d \mathbf v(t_1^{+})=\mathbf v(t_1^{+})-\mathbf v^*(t_1^{+*}) \notag 
\]
and $\delta(\cdot)$ is the variation, for example, $\delta \mathbf v_M(t_1^*)$ is the variation of $\mathbf v_M$ as an independent variable at $t_1^*$. 
The parts of $\delta \tilde J$ marked with asterisks are respectively due to the linear parts of the related integrals, and the parts of $\delta \tilde J$ marked with daggers are respectively obtained by the integrating by parts. 

Notice that by $\delta \tilde J=0$ and the fundamental lemma, it follows from the related integral terms in 
\eqref{eq_Dec5_01} that we first have the costate equation of the interceptor
\begin{align}
\dot{\mathbf p}_{Mri}(t)&=-\dfrac{\partial H_{M}(\mathbf r_M,\mathbf v_M, \mathbf p_{Mi})}{\partial \mathbf r_M}\notag\\\dot{\mathbf p}_{Mvi}(t)&=-\dfrac{\partial H_{Mi}(\mathbf r_M,\mathbf v_M, \mathbf p_{Mi})}{\partial \mathbf v_M}\notag
\end{align}
Due to the definition of Hamiltonian function \eqref{eq_7_2_1_Sep19}, we have
\begin{align}
\dot{\mathbf{p}}_{Mri}
&=-\dfrac{\mu}{r^3_M}(\dfrac{3}{r^2_M}\mathbf{r}_M\mathbf{r}^{\text T}_M-I_{3\times 3})\mathbf{p}_{Mvi}\notag\\
\dot{\mathbf{p}}_{Mvi}&=-\mathbf{p}_{Mri} \notag 
\end{align}
For the target, there are similar costate equations. We next derive all boundary conditions related to costates. By regrouping terms with respect to  $\delta \mathbf{v}_0, \delta \mathbf{v}_1$ in 
\eqref{eq_Dec5_01}, we have  
\[
\left(\dfrac{\Delta \mathbf{v}_0}{\abs{\Delta \mathbf{v}_0}}-\mathbf{q}_{v1}
\right)\delta \mathbf{v}_0, \quad 
\left(\dfrac{\Delta \mathbf{v}_1}{\abs{\Delta \mathbf{v}_1}}-\mathbf{q}_{v2}\right)
\delta \mathbf{v}_1
\]
Since  $\delta \mathbf{v}_0$ and $\delta \mathbf{v}_1$ are arbitrary,
their coefficients equal to zero. Hence we have BCs
 \begin{align}
\mathbf{q}_{v1}-\dfrac{\Delta \mathbf{v}_0}{\abs{\Delta \mathbf{v}_0}}=0, \quad 
 \mathbf{q}_{v2}-\dfrac{\Delta \mathbf{v}_1}{\abs{\Delta \mathbf{v}_1}} =0 \notag\end{align}
We make use of the relation between the difference $d(\cdot) $ and the variation $\delta(\cdot)$, for example, 
 \begin{align}
d \mathbf v(t_1^-)=\delta \mathbf v(t_1^{-*})+\dot{\mathbf v}(t_1^{-*})\delta t_1 \label{eq_Oct19_eq4} \end{align}
After a careful derivation,  by making the coefficients of the variations of all independent variables vanish remained in 
\eqref{eq_Dec5_01} such that $\delta \tilde J=0$, e.g., 
$
(\mathbf q_{v1}+\mathbf p_{v1}(t_0^+))d \mathbf{v}_0
$
where $d \mathbf{v}_0, d \mathbf{v}_M(t_0^+)$ have the same meaning,
we obtain the following BCs related to the costates
 \begin{align}
H_{M1*}(t_1^{-*})-H_{M2*}(t_1^{+*})+H_{T1*}(t_1^{-*})-H_{T2*}(t_1^{+*})=0\notag\\
H_{M2*}(t_h^{-*})+H_{T2*}(t_h^{-*})=0\notag\\
\mathbf q_{v1}+\mathbf p_{Mv1}(t_0^+)=0, \quad 
\mathbf q_{v2}+\mathbf p_{Mv1}(t_1^{-*})=0, \quad 
\mathbf q_{v2}+\mathbf p_{Mv2}(t_1^{+*})=0
\notag\\
\mathbf p_{Mv2}(t_h^{-*})=0\notag\\
\mathbf p_{Mr1}(t_1^{-*})-\mathbf p_{Mr2}(t_1^{+*})=0, \quad
\mathbf p_{Mr2}(t_h^{-*})-\bm\gamma_h=0 \notag\\
\mathbf p_{Tr1}(t_1^{-*})-\mathbf p_{Tr2}(t_1^{+*})=0, \quad 
\mathbf p_{Tr2}(t_h^{-*})+\bm\gamma_h=0\notag\\
\mathbf p_{Tv1}(t_1^{-*})-\mathbf p_{Tv2}(t_1^{+*})=0, \quad \mathbf p_{Tv2}(t_h^{-*})=0\notag
\end{align}
Notice that the continuity of the states at $t_h$ implies $d\mathbf r(t_h^{-})=d\mathbf r(t_h^+)=d\mathbf r(t_h)$ and so on. Furthermore, eliminating the intermediate quantities, e.g., $\bm\gamma_h$, we obtain BCs
\begin{align}
&\begin{cases}
\mathbf p_{Mv1}(t_0^+)+\dfrac{\Delta \mathbf{v}_0}{\abs{\Delta \mathbf{v}_0}}=0, \quad \mathbf p_{Mv1}(t_1^{-*})+
\dfrac{\Delta \mathbf{v}_1}{\abs{\Delta \mathbf{v}_1}} =0\notag\\
\mathbf p_{Mv1}(t_1^{-*})-\mathbf p_{Mv2}(t_1^{+*})=0, \quad
\mathbf p_{Mv2}(t_h^{-*})=0 \notag\\
\end{cases}\\
&\begin{cases}
H_{M1*}(t_1^{-*})-H_{M2*}(t_1^{+*})+H_{T1*}(t_1^{-*})-H_{T2*}(t_1^{+*})=0\notag\\
H_{M2*}(t_h^{-*})+H_{T2*}(t_h^{-*})=0\notag\\
\end{cases}\\
&\mathbf p_{Mr1}(t_1^{-*})-\mathbf p_{Mr2}(t_1^{+*})=0, \quad
\mathbf p_{Mr2}(t_h^{-*})+
\mathbf p_{Tr2}(t_h^{-*})=0
\notag\\
&\mathbf p_{Tr1}(t_1^{-*})-\mathbf p_{Tr2}(t_1^{+*})=0, \quad 
\mathbf p_{Tv1}(t_1^{-*})-\mathbf p_{Tv2}^{\text T}(t_1^{+*})=0, \quad \mathbf p_{Tv2}(t_h^{-*})=0\notag
\end{align}
In order to obtain explicit BCs, we now calculate the involved Hamiltonian functions. 
Since
\begin{align}
H_{T2*}(t_h^{*})&=\mathbf p_{Tr2}^{\text T}(t_h^*)\mathbf v_T(t_h^*)-\dfrac{\mu}{r_T^2(t_h^*)}\mathbf p_{Tv2}^{\text T}(t_h^*)\mathbf r_{T}(t_h^*)\notag\\
H_{M2*}(t_h^{*})&=\mathbf p_{Mr2}^{\text T}(t_h^*)\mathbf v_M(t_h^*)-\dfrac{\mu}{r_M^2(t_h^*)}\mathbf p_{Mv2}^{\text T}(t_h^*)\mathbf r_{M}(t_h^*)\notag
\end{align}
we have
\begin{align}
H_{M2*}(t_h^{*})+
H_{T2*}(t_h^{*})
&=\mathbf p_{Tr2}^{\text T}(t_h^*)\mathbf v_T(t_h^*)
+\mathbf p_{Mr2}^{\text T}(t_h^*)\mathbf v_M(t_h^*)\notag\\
&=\mathbf p_{Tr2}^{\text T}(t_h^*)\left(\mathbf v_T(t_h^*)
-\mathbf v_M(t_h^*)\right)=\mathbf p_{Mr2}^{\text T}(t_h^*)\left(\mathbf v_M(t_h^*)-
\mathbf v_T(t_h^*)\right)=0\notag
\end{align}
It follows from the continuity of the target state that
\begin{align}
H_{T1*}(t_1^{-*})-H_{T2*}(t_1^{+*})=0\notag
\end{align}
Hence 
\begin{align}
H_{M1*}(t_1^{-*})-H_{M2*}(t_1^{+*})+H_{T1*}(t_1^{-*})-H_{T2*}(t_1^{+*})&=H_{M1*}(t_1^{-*})-H_{M2*}(t_1^{+*})\notag\\
&=\mathbf p_{Mr1}^{\text T}(t_1^{-*})\mathbf v_M(t_1^{-*})-
\mathbf p_{Mr2}^{\text T}(t_1^{+*})\mathbf v_M(t_1^{+*})\notag\\
&=-\mathbf p_{Mr1}^{\text T}(t_1^{-*})\Delta \mathbf v(t_1)\notag\\
&=-\mathbf p_{Mr1}^{\text T}(t_1^{+*})\Delta \mathbf v(t_1)
=0\notag
\end{align}
We finally obtain BCs for Problem \ref{problem1} related to costates
\begin{align}
&(1)\begin{cases}
\mathbf p_{Mv1}(t_0^+)+\dfrac{\Delta \mathbf{v}_0}{\abs{\Delta \mathbf{v}_0}}=0, \quad \mathbf p_{Mv1}(t_1^{-*})+
\dfrac{\Delta \mathbf{v}_1}{\abs{\Delta \mathbf{v}_1}} =0\notag\\
\mathbf p_{Mv1}(t_1^{-*})-\mathbf p_{Mv2}(t_1^{+*})=0, \quad
\mathbf p_{Mv2}(t_h^{-*})=0 \notag\\
\mathbf p_{Mr1}(t_1^{-*})-\mathbf p_{Mr2}(t_1^{+*})=0\notag
\end{cases}\\
&(2)\begin{cases}
-\mathbf p_{Mr1}^{\text T}(t_1^{-*})\Delta \mathbf v(t_1)=0
\notag\\
\mathbf p_{Tr2}^{\text T}(t_h^*)\left(\mathbf v_T(t_h^*)
-\mathbf v_M(t_h^*)\right)=0 
\notag\\
\end{cases}\\
&(3)~\mathbf p_{Mr2}(t_h^{-*})+
\mathbf p_{Tr2}(t_h^{-*})=0
\notag\\
&(4) \quad \mathbf p_{Tr1}(t_1^{-*})-\mathbf p_{Tr2}(t_1^{+*})=0, \quad 
\mathbf p_{Tv1}(t_1^{-*})-\mathbf p_{Tv2}^{\text T}(t_1^{+*})=0, \quad \mathbf p_{Tv2}(t_h^{-*})=0\notag
\end{align}
The second group is from the Hamiltonian conditions. 
Actually in numerical experiments, the boundary conditions (3)-(4) related to the costates of the target may not be used if we do not consider the dynamic of the target which is just included in BCs to provide its position and velocity.  
Notice that in the sequel, to simplify notations, the ``$*$'' symbol  which denotes an optimal value has been removed. The boundary conditions related to the costates of the interceptor is shown in List \ref{table_Nov7_01}. 

\section*{Appendix B: Initial Data} \label{App_2}

The examples in this paper involve three groups of the initial data of the interceptor and target generated by ballistic missiles, which are numerically obtained by corresponding orbital elements. 
The orbital elements  of the interceptor and the target are shown in Table \ref{table04c}. 
\begin{center}
\begin{threeparttable}\small 
\caption{Orbital elements}\label{table04c}
\begin{tabular}{p{2cm}c ccccc}
\hline
& spacecraft & $H$ & $i$ & $\Omega$ & $e$ & $\omega $ \\
\hline
first group & interceptor & 500 km & $2\pi/3$ & $4\pi/3$ & 0.3 & $\pi$ \\
& target  & 400 km & $\pi/6$ &  $\pi/3$ & 0.1 & 0 
\\ \hline
second group & interceptor & 500 km & $2\pi/3$ & $5\pi/4$ & 0.5 & $\pi$ \\
& target & 400 km & $\pi/6$ &  $\pi/3$ & 0.2 & 0 
\\ \hline
\end{tabular}
\end{threeparttable}
\end{center}

For the case the true anomaly $\theta=0$,  by using Algorithm 4.5 in \cite{Curits_book}, 
the trajectories of the interceptor and the target are known, and we then obtain  
the initial positions of the interceptor and the target corresponding to the first time at which the target and the interceptor are 
above the atmosphere (approximately equal to the sum of the earth's radius and the height of the atmosphere 120 km) and their velocities. For the first group of the orbital elements, we generate the initial data I and II, and for the second group, we have the initial data III. A reference vector of the interceptor is chosen just above the atmosphere for the initial data I and II.
{\small 
\begin{verbatim}
The initial data I:
the initial position of the target
   1.0e+06*[-5.842891129580837; -1.241946037180446; 2.562926625347858]
the initial velocity of the target
   1.0e+03*[-0.065508668182581; -7.322759468283627; -2.081144241020925]
the initial position of the interceptor
    1.0e+06*[-1.392985266715916; -5.682521353135304; -2.831729949288823]
the initial velocity of the interceptor
    1.0e+03*[-4.511678481085538; -2.680368719222989; 4.446250319272038]
the reference position of the interceptor
    1.0e+06*[-4.4528; -4.4166; 1.7258]
--------------
The initial data II:
the initial position of the target
   1.0e+06*[-5.842481237484495; -1.389922138771051; 2.520004658256203]
the initial velocity of the target
   1.0e+03*[0.105801179312784; -7.284177899593129; -2.155661625234000]
the initial position of the interceptor
    1.0e+06*[-1.422033750436706; -5.699632649250217; -2.802976040834825]
the initial velocity of the interceptor
    1.0e+03*[-4.498532928342012; -2.627216111719469; 4.472563498532185]
--------------
The initial data III:
the initial position of the target
   1.0e+06*[-5.394452557207117; -3.192217335202957; 1.775712509707950]
the initial velocity of the target
   1.0e+03*[1.767918629073472; -6.417485911429783; -2.736527923779045]
the initial position of the interceptor
   1.0e+06*[-3.580084601432768; -5.106010405266481; -1.868869802369432]
the initial velocity of the interceptor
   1.0e+03*[-4.211400455599469; -0.835510934866194; 4.134603376902692]
\end{verbatim}
}

\section*{Appendix C: BCs, parameters, and initial values for Examples \ref{exam1}-\ref{exam2}} \label{App_31}

\begin{bcs} \label{table_Oct10_01}
\medskip
\begin{center} 
\centering{\rm \small BCs related to the costates and inequality constraints\vspace*{0.2cm}}

{\blue
\fbox{
\begin{minipage}[H]{10cm} \vspace*{-.3cm}
\black\small
\begin{align}
\begin{array}{ll}
&(1)\begin{cases}
\mathbf p_{Mv1}(t_1^-)+\dfrac{\Delta \mathbf{v}_1}{\abs{\Delta \mathbf{v}_1}}+
\begin{bmatrix} \mu_1-\mu_2 & \mu_3-\mu_4 &  \mu_5-\mu_6\end{bmatrix}^{\text T}
=0\notag\\
\mathbf p_{Mv2}(t_2^-)+
\dfrac{\Delta \mathbf{v}_2}{\abs{\Delta \mathbf{v}_2}} 
+
\begin{bmatrix} \mu_7-\mu_8 & \mu_9-\mu_{10} &  \mu_{11}-\mu_{12}\end{bmatrix}^{\text T}
=0\notag\\
\mathbf p_{Mv1}(t_1^-)-\mathbf p_{Mv2}(t_1^+)=0, \quad
\mathbf p_{Mv2}(t_2^{-*})-\mathbf p_{Mv3}(t_2^{+*})=0 
\notag\\
\mathbf p_{Mr1}(t_1^{-*})-\mathbf p_{Mr2}(t_1^{+*})=0, \quad 
\mathbf p_{Mr2}^{\text T}
(t_2^{-*})-\mathbf p_{Mr3}^{\text T}(t_2^{+*})=0\notag\\
\mathbf p_{Mv3}(t_h^{-*})-\mathbf p_{Mv4}(t_h^{+*})=0
\notag
\end{cases}\vspace*{0.1cm}\\
&(2)\begin{cases}
-\mathbf p_{Mr1}^{\text T}(t_1^{-*})\Delta \mathbf v(t_1)-\lambda_1+\lambda_2+\lambda_3=0\notag\\
-\mathbf p_{Mr2}^{\text T}(t_2^{-*})\Delta \mathbf v(t_2)-\lambda_3=0\notag\\
\mathbf p^{\text T}_{Mr3}(t_h^{-*})
\left(\mathbf v_M(t_h^*)-\mathbf v_T(t_h^*)\right)=0
\notag
\end{cases}\vspace*{0.1cm}\\
&(3)\begin{cases}\vspace*{0.1cm}
\diag(\mu_1, \mu_3, \mu_5)
\left(\Delta \mathbf{v}_1-\mathbf p_{1max}\right)=0
\\\vspace*{0.1cm}
\diag(\mu_2,\mu_4,\mu_6)
\left(\mathbf p_{1min}
-\Delta \mathbf v_1\right)=0\\\vspace*{0.1cm}
\diag(\mu_7,\mu_9,\mu_{11})
\left(\Delta \mathbf{v}_2-\mathbf p_{2max}\right)=0
\\\vspace*{0.1cm}
\diag(\mu_8, \mu_{10}, \mu_{12})
\left(\mathbf p_{2min}
-\Delta \mathbf v_2\right)=0
\end{cases}\vspace*{0.1cm}
\end{array}
\notag 
\end{align}
\end{minipage}}}
\end{center}
\end{bcs}

\medskip
The parameters and initial values of Examples \ref{exam1}-\ref{exam2} are given in Tables \ref{table06b}-\ref{table06a}. 
\begin{center}
\begin{threeparttable}\small 
\caption{Parameters and initial values of Example \ref{exam1}}\label{table06b}
\begin{tabular}{c cccc  c ccccc }
\toprule
\multirow{2}{*}{parameters}& \multirow{2}{*}{$\alpha$} &\multirow{2}{*}{$\beta$} &\multirow{2}{*}{$\gamma$} & \multirow{2}{*}{tolerance} & $ p_{1min}$ & $ p_{2min}$  & $ p_{3min}$  & $ p_{4min}$  & $ p_{5min}$  & $ p_{6min}$  \\
&   &  &  & & $ p_{1max}$ & $ p_{2max}$  & $ p_{3max}$  & $ p_{4max}$  & $ p_{5max}$  & $ p_{6max}$  \\
\hline
&  \multirow{2}{*}{20} &  \multirow{2}{*}{40} &  \multirow{2}{*}{50} &  \multirow{2}{*}{1e-12} & -1300 &-1300& -1300 & -1300 & -1300 &-1300 \\
&  & & &  &  1300 & 1300 & 1300 & 1300 & 1300 &1300
\\ \hline
initial values  &  $t_1$ &  {$t_2$} &  {$t_h$} &  {$\lambda_1, \lambda_2, \lambda_3$} &  \multicolumn{6}{c} {$\mu_1, \mu_2, \dots,\mu_{12}$}  \\
\hline
& 0.04   & 0.14 & 500 & 1 &  \multicolumn{6}{c}{1}
\\ 
\bottomrule
\end{tabular}
\end{threeparttable}
\end{center}
\begin{center}
\begin{threeparttable}\small 
\caption{Parameters and initial values of Example \ref{exam2}}\label{table06a}
\begin{tabular}{c cccc c ccccc }
\toprule
\multirow{2}{*}{parameters}& \multirow{2}{*}{$\alpha$} &\multirow{2}{*}{$\beta$} &\multirow{2}{*}{$\gamma$} & \multirow{2}{*}{tolerance} & $ p_{1min}$ & $ p_{2min}$  & $ p_{3min}$  & $ p_{4min}$  & $ p_{5min}$  & $ p_{6min}$  \\
&   &  &  & & $ p_{1max}$ & $ p_{2max}$  & $ p_{3max}$  & $ p_{4max}$  & $ p_{5max}$  & $ p_{6max}$  \\
\hline
&  \multirow{2}{*}{20} &  \multirow{2}{*}{40} &  \multirow{2}{*}{50} &  \multirow{2}{*}{1e-12} & -400 &-400& -500 & -400 & -400 &-500 \\
&  & & &  &  400 & 400 & 400 & 400 & 400 &400
\\ \hline
initial values  &  $t_1$ &  {$t_2$} &  {$t_h$} &  {$\lambda_1, \lambda_2, \lambda_3$} & \multicolumn{6}{c} {$\mu_1, \mu_2, \dots,\mu_{12}$} \\
\hline
& 0.03   & 0.1 & 700 & 1 & \multicolumn{6}{c}{1}
\\ 
\bottomrule
\end{tabular}
\end{threeparttable}
\end{center}
\medskip

In both examples, we provide the initial guesses for velocity impulses and costates
{\small
\begin{verbatim}
dv10=[-300; 300; -300], dv20=[-100; 100; -100]
pmr0=1.0e-03*[0.5396; -0.4740; 0.8451], pmv0=[0.4862; -0.4365; 0.7570]
\end{verbatim}
}

\section*{Appendix D: Time Change} \label{App_3}

We have transformed two-impulse space interception problems into multi-point boundary value problems by using the calculus of variations. 
In order to use the Matlab boundary value problem solvers bvp4c and bvp5c to solve them, we must introduce those instants $t_1,t_2,t_h,t_f$ as unknown parameters because the solvers cannot directly deal with a multi-point BC value problem with unknown boundary points of time intervals. 
The solvers require the normalization of the unknown final time and the parameterization of unknown time instants. For example, consider the two-impulse space interception problem \ref{problem1}. 
We introduce a time change
\begin{align}
\tau=\dfrac{t-t_0}{t_h-t_0}, \quad t\in [t_0, t_h] \label{eq_Sep20_1}
\end{align}
then $t_1, t_2,t_h$ are transformed into 
\begin{align}
\tau_1=\dfrac{t_1-t_0}{t_h-t_0}, \quad \tau_2=\dfrac{t_2-t_0}{t_h-t_0}, \quad \tau_h=1
\end{align}
In this way, the time change \eqref{eq_Sep20_1} transforms $[0, t_1]$ into $[0, \tau_1]$ and so on.  The non-dimensional factors $\tau_i$ are call scaled time instants. 
We now introduce time changes for each sub-intervals
\begin{align}
s=\begin{cases}
\dfrac{\tau}{3\tau_1}, \quad &\tau\in [0, \tau_1] \vspace*{0.2cm}\\
\dfrac{\tau-2\tau_1+\tau_2}{3(\tau_2-\tau_1)}, \quad &\tau\in [\tau_1, \tau_2] \vspace*{0.2cm}\\
\dfrac{\tau-3\tau_2+2}{3(1-\tau_2)}, \quad & \tau\in [\tau_2, 1] 
\end{cases}\label{eq_Sep20_02}
\end{align}

Under the time $t$, denote $x=\begin{bmatrix}\mathbf r & \mathbf v\end{bmatrix}^\text{T}$, $f(x)=\begin{bmatrix}{\mathbf v} & -\dfrac{\mu}{r^3}\mathbf r\end{bmatrix}^\text{T}$, we have 
\begin{align}\label{eq_Dec26_01}
\dfrac{d x}{dt}=\begin{cases}
f(x), \quad t\in[t_0, t_1)\\
f(x), \quad t\in[t_1, t_2)\\ 
f(x), \quad t\in[t_2, t_h]
\end{cases}
\end{align}
With the time changes \eqref{eq_Sep20_1}-\eqref{eq_Sep20_02}, we have 
\begin{align}
\dfrac{d x}{dt}=\begin{cases}
f(x), \quad t\in[t_0, t_1)\\
f(x), \quad t\in[t_1, t_2)\\ 
f(x), \quad t\in[t_2, t_h]
\end{cases}
\Longrightarrow 
\dfrac{d x}{ds}=\begin{cases}
3t_h\tau_1f(x), \quad & s\in[0, 1/3)\\
3t_h(\tau_2-\tau_1)f(x), \quad & s\in[1/3, 2/3)\\ 
3t_h(1-\tau_2)f(x), \quad & s\in[2/3, 1]
\end{cases}
\end{align}
After introducing the time changes, the boundary conditions at $t_1$ are transformed into the boundary conditions at $1/3$ and so on. 
Similar time changes can be applied to other problems, e.g.~Problem \ref{problem2} in which there are four sub-intervals.
Then the solvers can be used to solve the multi-point boundary value problem of the above piecewise continuous ODEs with unknown parameters. 
The number of BCs is the sum of  the number of sub-intervals multiplied by the number of equations and the number of unknown parameters. For time changes,  we refer to, e.g., \cite[Appendix A]{Longuski_book} and \cite[Section 3]{Zefran} for details.    

\section*{Acknowledgments}
The first author is supported by the National Natural Science Foundation of China (no. 61374084). We appreciate the help of Professor Prussing and Dr.~Sandrik for sharing their Matlab scripts in \cite{Sandrik2006}. Thanks also go to Dr.~Kierzenka, one of the authors of these  Matlab boundary value problem solvers, for sending us his Ph.D.~Dissertation \cite{Kierzenka_thesis}.  The reviewers' comments help us improve the presentation of the manuscript. We appreciate the Matlab solvers, and without them we do not have this work.

\bibliographystyle{plain}

\begin{thebibliography}{27}
\newcommand{\enquote}[1]{``#1''}
\providecommand{\natexlab}[1]{#1}
\providecommand{\url}[1]{\texttt{#1}}
\providecommand{\urlprefix}{URL }
\expandafter\ifx\csname urlstyle\endcsname\relax
  \providecommand{\doi}[1]{doi:\discretionary{}{}{}#1}\else
  \providecommand{\doi}{doi:\discretionary{}{}{}\begingroup
  \urlstyle{rm}\Url}\fi

\bibitem[{{This Week's Citation Classic}(1980)}]{Bryson_story}
{This Week's Citation Classic}, February 25 1980.
\newblock\\
  \urlprefix\url{http://garfield.library.upenn.edu/classics1980/A1980JE96000001.pdf}.

\bibitem[{{Bryson Jr}(1996)}]{Bryson1996}
{Bryson Jr}, A.~E., \enquote{Optimal Control--1950 to 1985,} \emph{IEEE Control
  Systems}, 1996, pp. 26--33.

\bibitem[{Tsien and Evans(1951)}]{Tsien}
Tsien, H.~S., and Evans, R.~C., \enquote{Optimum Thrust Programming for a
  Sounding Rocket,} \emph{Journal of the American Rocket}, Vol.~21, No.~5,
  1951, pp. 99--107.

\bibitem[{Lawden(1963)}]{Lawden_book}
Lawden, D.~F., \emph{Optimal Trajectories for Space Navigation}, Butter Worths,
  London, 1963.

\bibitem[{Bryson~Jr. and Ho(1975)}]{Bryson_book}
Bryson~Jr., A.~E., and Ho, Y.-C., \emph{Applied Optimal Control}, Hemisphere
  Publishing Corp., London, 1975.

\bibitem[{James M.~Longuski and Prussing(2014)}]{Longuski_book}
James M.~Longuski, J. J.~G., and Prussing, J.~E., \emph{Optimal Control with
  Aerospace Applications}, Springer, 2014.

\bibitem[{Ben-Asher(2010)}]{Ben_Asher_book}
Ben-Asher, J.~Z., \emph{Optimal Control Theory with Aerospace Applications},
  AIAA, 2010.

\bibitem[{Subchan and \.{Z}bikowski(2009)}]{Subchan_book}
Subchan, S., and \.{Z}bikowski, R., \emph{Computational Optimal Control},
  Wiley, 2009.

\bibitem[{Curtis(2014)}]{Curits_book}
Curtis, H.~D., \emph{Orbital Mechanics for Engineering Students},
  3\textsuperscript{rd} ed., Elsevier Ltd., 2014.

\bibitem[{Prussing(2010)}]{prussing_2010}
Prussing, J.~E., \enquote{Primer Vector Theory and Applications,}
  \emph{Spacecraft Trajectory Optimization}, edited by B.~A. Conway, Cambridge,
  2010, Chap.~2, pp. 16--36.

\bibitem[{Gobetz and Doll(1969)}]{Gobetz1969}
Gobetz, F.~W., and Doll, J.~R., \enquote{A Survey of Impulsive Trajectories,}
  \emph{AIAA Journal}, Vol.~7, No.~5, 1969, pp. 801--834.

\bibitem[{Robinson(1967)}]{Robinson_1967}
Robinson, A.~C., \enquote{Comparison of Fuel-optimal Maneuvers Using a Minimum
  Number of Impulses with Those Using the Optimal Number of Impulses-A Survey,}
  {NASA} contractor report {NASW}-1146, Battelle Memorial Inst., 1967.

\bibitem[{Jezewski(1975)}]{Jezewski1975}
Jezewski, D.~J., \enquote{Primer Vector Theory and Applications,} Tech. rep.,
  NASA TR R-454, November 1975.

\bibitem[{Vinh et~al.(1990)Vinh, Lu, Howe, and Gilbert}]{Vinh1990}
Vinh, N.~X., Lu, P., Howe, R.~M., and Gilbert, E.~G., \enquote{Optimal
  Interception with Time Constraint,} \emph{Journal of Optimization Theory and
  Applications}, Vol.~66, No.~3, 1990, pp. 361--390.

\bibitem[{Taur et~al.(1995)Taur, Coverstone-Carroll, and Prussing}]{Taur1995}
Taur, D.-R., Coverstone-Carroll, V., and Prussing, J.~E., \enquote{Optimal
  Impulsive Time-Fixed Orbital Rendezvous and Interception with Path
  Constraints,} \emph{Journal of Guidance, Control, and Dynamics}, Vol.~18,
  No.~1, 1995, pp. 54--60.

\bibitem[{Prussing and Chiu(1986)}]{Prussing1986}
Prussing, J.~E., and Chiu, J.-H., \enquote{Optimal Multiple-Impulse Time-Fixed
  Rendezvous between Circular Orbits,} \emph{Journal of Guidance, Control, and
  Dynamics}, Vol.~9, No.~1, 1986, pp. 17--22.

\bibitem[{Sandrik(2006)}]{Sandrik2006}
Sandrik, S., \enquote{Primer-Optimized Results and Trends for Circular Phasing
  and Other Circle-to-Circle Impulsive Coplanar Rendezvous,} Ph.D. thesis,
  University of Illinois at Urbana-Champaign, 2006.

\bibitem[{Colasurdo and Pastrone(1994)}]{Colasurdo1994}
Colasurdo, G., and Pastrone, D., \enquote{Indirect Optimization Method for
  Impulsive Transfers,} \emph{Astrodynamics Conference, Guidance, Navigation,
  and Control}, 1994, pp. 441--448.

\bibitem[{Prussing(1995)}]{Prussing1995}
Prussing, J.~E., \enquote{Optimal Impulsive Linear Systems: Sufficient
  Conditions and Maximum Number of Impulses,} \emph{The Journal of the
  Astronautical Sciences}, Vol.~43, No.~2, 1995, pp. 195--206.

\bibitem[{Shampine et~al.(2003)Shampine, Gladwell, and
  Thompson}]{Shampine_book}
Shampine, L.~F., Gladwell, I., and Thompson, S., \emph{Solving {ODE}s with
  {M}atlab}, Cambridge University Press, 2003.

\bibitem[{Kierzenka(1998)}]{Kierzenka_thesis}
Kierzenka, J., \enquote{Studies in the Numerical Solution of Ordinary
  Differential Equations,} Ph.D. thesis, Department of Mathematics, Southern
  Methodist University, Dallas, TX., 1998.

\bibitem[{Kierzenka and Shampine(2008)}]{Kierzenka_2008}
Kierzenka, J.~A., and Shampine, L.~F., \enquote{A BVP Solver that Controls
  Residual and Error,} \emph{JNAIAM J. Numer. Anal. Ind. Appl. Math}, Vol.~3,
  2008, pp. 27--41.

\bibitem[{Kierzenka and Shampine(2001)}]{Kierzenka_2001}
Kierzenka, J.~A., and Shampine, L.~F., \enquote{A BVP Solver Based on Residual
  Control and the MATLAB PSE,} \emph{ACM Transactions on Mathematical
  Software}, Vol.~27, No.~3, 2001, pp. 299--316.

\bibitem[{Sigal and Ben-Asher(2014)}]{Ben_Asher2014}
Sigal, E., and Ben-Asher, J.~Z., \enquote{Optimal Control for Switched Systems
  with Pre-defined Order and Switch-Dependent Dynamics,} \emph{Journal of
  Optimization Theory and Applications}, Vol. 161, No.~2, 2014, pp. 582--591.

\bibitem[{Xie et~al.(2017)Xie, Zhang, and Xu}]{JOTA_XZX}
Xie, L., Zhang, Y., and Xu, J., \enquote{Hohmann Transfer via Constrained
  Optimization,} {https:}//arxiv.org/abs/1712.01512, December 2017.

\bibitem[{Hull(2003)}]{Hull_book}
Hull, D.~G., \emph{Optimal Control Theory for Applications}, Springer, 2003.

\bibitem[{Zefran et~al.(1996)Zefran, Desai, and Kumar}]{Zefran}
Zefran, M., Desai, J.~P., and Kumar, V., \enquote{Continuous Motion Plans for
  Robotic Systems with Changing Dynamic Behavior,} \emph{The second Int.
  Workshop on Algorithmic Foundations of Robotics}, Toulouse, France, 1996.

\end{thebibliography}

\end{document}